\documentclass[a4paper,english,oneside]{article}
\usepackage[T1]{fontenc}
\usepackage[utf8]{inputenc} 

\usepackage{geometry}
\usepackage{amsmath}
\usepackage{amsfonts}
\usepackage{amssymb}
\usepackage{amsthm}
\usepackage[english]{babel}
\usepackage{color}
\usepackage{lmodern}
\usepackage{pstricks-add}
\usepackage{dsfont}
\usepackage{hyperref}
\hypersetup{hidelinks}

\usepackage{mathscinet}
\usepackage{enumitem}
\usepackage[cal=boondoxo,bb=ams]{mathalfa}
\usepackage{microtype}

\usepackage{mathtools}
\usepackage{esint}

\numberwithin{equation}{section}

\theoremstyle{plain}
\newtheorem{theorem}{Theorem}[section]

\theoremstyle{definition}
\newtheorem{definition}[theorem]{Definition}

\theoremstyle{remark}
\newtheorem{remark}{Remark}

    \DeclareMathOperator\supp{supp}

    \DeclareMathOperator\meas{meas}
 
 \date{}
     
\begin{document}

\title{A note on blow-up results for semilinear wave equations in de Sitter and anti-de Sitter spacetimes}

\author{Alessandro Palmieri$\,^{\mathrm{a}}$, Hiroyuki Takamura$\,^{\mathrm{b}}$}

\date{
\small{ $\,^\mathrm{a}$ Mathematical Institute, Tohoku University, Aoba, Sendai 980-8578, Japan} \\
\small{ $\,^\mathrm{b}$ Mathematical Institute/Research Alliance Center of Mathematical Sciences, Tohoku University, Aoba, Sendai 980-8578, Japan}  }

\maketitle
\begin{center}
\normalsize{\today}
\end{center}

\begin{abstract}
In this work we derive some blow-up results for semilinear wave equations both in de Sitter and anti-de Sitter spacetimes. By requiring suitable conditions on a time-dependent factor in the nonlinear term, we prove the blow-up in finite time of the spatial averages of local in time solutions. In particular, we derive a sequence of lower bound estimates for the spatial average by combining a suitable slicing procedure with an iteration frame for this time-dependent functional.
\end{abstract}

\begin{flushleft}
\end{flushleft}

\begin{flushleft}
\textbf{AMS Classification (2020)}  Primary: 35L05, 35L71, 35B44; Secondary:  33C10.
\end{flushleft}

\section{Introduction}

In the first part of this paper, we investigate the blow-up dynamic for local in time solutions to the following semilinear wave equation with damping and mass in de Sitter spacetime
\begin{align}\label{deSitter semi} 
\begin{cases}
\partial_t^2 u- c^2\mathrm{e}^{-2Ht} \Delta u+ b\partial_t u +m^2 u=f(t,u), & x\in \mathbb{R}^n, \ t\in (0,T), \\
u(0,x)= \varepsilon u_0(x), & x\in \mathbb{R}^n, \\
\partial_t u(0,x)= \varepsilon u_1(x), & x\in \mathbb{R}^n,
\end{cases}
\end{align}  where $c,H$ are positive constants, $b,m^2$ are nonnegative real parameters, $\varepsilon>0$ is a parameter describing the size of initial data and $T=T(\varepsilon)\in (0,\infty]$ is the lifespan (maximal existence time) of the weak solution with $C^1$ regularity in time variable (cf. Definition \ref{Def sol}). In the literature concerning cosmology, the constant $H$ is called Hubble constant, $m$ denotes the mass of a particle and the coefficient $b$ for the damping term is taken equal to the space dimension $n$ (see, for example, (0.6) in \cite{YagGal09}).

As nonlinear term we consider a nonlocal term given by the product of three terms: a time-dependent coefficient providing (possible) additional exponential and/or polynomial growth, a $p$-power nonlinearity, and a power of the spatial $L^p$ norm of the solution. Namely, we set
\begin{align}\label{nonlinearity}
f(t,u)\doteq \Gamma(t) \left(\int_{\mathbb{R}^n}|u(t,y)|^p\mathrm{d}y\right)^\beta |u|^p,
\end{align} where $p>1$, $\beta\geqslant 0$, and $\Gamma(t)$ is a suitable nonnegative function. Our goal in the present paper will be to determine growth conditions on $\Gamma=\Gamma(t)$ (depending on $p,\beta,b,m^2$) in a such way that blow-up phenomena for the local solutions to \eqref{deSitter semi} occur under suitable sign assumptions for the Cauchy data.

We investigate the case in which the damping term $b\partial_t u$ is dominant over the mass term $m^2u$, by prescribing a restriction on the size of $m^2$. More precisely, we will always work under the following assumption
\begin{align} \label{dominant damping}
b^2\geqslant 4m^2
\end{align} for the coefficients of the lower order terms. Following the nomenclature introduced in \cite{ER18}, we call $b^2 > 4m^2$ the case with \emph{dominant dissipation}, $b^2 = 4m^2$ the case with \emph{balanced dissipation and mass}, and $b^2 < 4m^2$ the case with \emph{dominant mass}. We do not consider the dominant mass case since this case is somehow related to Klein-Gordon equation with real positive mass, which cannot be treated with the approach that we are going to use in the present work. More specifically, we investigate the dynamic of the space average of a local solution to \eqref{deSitter semi}, by determining a lower bound estimate for this time-dependent functional, where the space average appears also in a nonlinear form in an integral term on the right-hand side of this inequality (the so-called \emph{iteration frame}). For this kind of approach it is essential to work with nonnegative lower bounds and \eqref{dominant damping} ensures us that the time-dependent factors on the right-hand side of the iteration frame have no oscillations and are positive.

We point out that the speed of propagation, namely, the function $a_{\mathrm{dS}}(t)\doteq c\, \mathrm{e}^{-Ht}$, is exponentially decreasing in the previous semilinear wave equation. Moreover, the amplitude of the forward light-cone, provided by $$A_{\mathrm{dS}}(t)\doteq \int_0^t a_{\mathrm{dS}}(\tau) \mathrm{d}\tau=\frac{c}{H}\left(1-\mathrm{e}^{-Ht}\right)$$ is a bounded function. 
In other words, by working with smooth solutions, if we assume $u_0$ and $u_1$ compactly supported in $B_R\doteq \{x\in\mathbb{R}^n:|x|\leqslant R\}$, given a local solution $u$ to \eqref{deSitter semi}, we have that
\begin{align}\label{support condition sol}
\supp u(t,\cdot) \subset B_{R+A_{\mathrm{dS}}(t)} \ \ \mbox{for any} \ t\in (0,T).
\end{align} For this support condition we used the property of finite speed of propagation or, alternatively, the explicit representation formulas from the series of works by Galstian and Yagdjian \cite{YagGal08,YagGal09,Yag09,Yag10,Yag13,Yag15}.
Therefore, assuming compactly supported Cauchy data, the support of a local in time solution will be contained in an infinite half cylinder (as long as the solution exists). As we will see in the proof of our blow-up results, this property will play a key role when establishing the iteration frame. 

We emphasize that the inclusion of the time-dependent factor $\Gamma$ in \eqref{nonlinearity} is made in order to be able to prove the blow-up in finite time for $p>1$ and $\beta\geqslant 0$. Indeed, the exponentially decaying speed of propagation and the presence of the mass term both make extremely difficult the occurrence of a blow-up in finite time of the solution. In the massless case (i.e. for $m^2=0$), we will be able not to require any additional exponential growth in the nonlinear term, that is, we may consider the case $\Gamma(t)=1$ as well. In particular, for $\beta=0$ we will recover (with a different technique) the result recently proved by Tsutaya-Wakasugi \cite{TsuWa21} with the test function method. On the contrary, when a mass term is present in the partial differential operator on the left-hand side of \eqref{deSitter semi} and we work under the assumption \eqref{dominant damping}, then, our method produces a sequence of lower bound estimates too weak that is not enough to prove the blow-up in finite time unless we require additional exponential growth through the factor $\Gamma(t)$.

We emphasize that the local case $\beta=0$ can be included in our result as well. This corresponds to the usual power nonlinearity with a time-dependent factor.

The nonlinear term in \eqref{nonlinearity} has been already considered in the literature for the Klein-Gordon equation in de Sitter spacetime by Yagdjian \cite{Yag09} when \eqref{dominant damping} is satisfied and by Nakamura \cite{Na20} for a pure imaginary mass (i.e. for $m^2<0$ with our notations) both for de Sitter and anti-de Sitter spacetimes. In both these papers a blow-up result is proved by means of a comparison argument for a certain ODE. In our approach we work with the corresponding integral formulation that will allow us to slightly improve the growth condition for $\Gamma(t)$ in comparison to that one in \cite[Theorem 1.1]{Yag09}. 

In the second part of the paper, we investigate what happens if we consider an exponentially increasing speed of propagation, say $a_{\mathrm{AdS}}(t)=c\, \mathrm{e}^{Ht}$ with $c,H>0$, in place of an exponentially decreasing function as in \eqref{deSitter semi}. In other words, we are interested to study the following semilinear problem associated with the \emph{wave equation in anti-de Sitter spacetime}
\begin{align}\label{anti deSitter semi} 
\begin{cases}
\partial_t^2 v- c^2\mathrm{e}^{2Ht} \Delta v+ b\partial_t v +m^2 v=f(t,v), & x\in \mathbb{R}^n, \ t\in (0,T), \\
v(0,x)= \varepsilon v_0(x), & x\in \mathbb{R}^n, \\
\partial_t v(0,x)= \varepsilon v_1(x), & x\in \mathbb{R}^n,
\end{cases}
\end{align}  where $c,H$ are positive constants, $b,m^2$ are nonnegative real parameters satisfying \eqref{dominant damping} and the nonlinear term is defined analogously as in  \eqref{nonlinearity}.

As in the corresponding results for the semilinear wave equation in de Sitter spacetime, we want to examine the growth conditions on the factor $\Gamma$ that provide local in time solutions that blow up in finite time (under suitable sign conditions for the Cauchy data). To the knowledge of the authors, while \eqref{deSitter semi} has been already studied in the literature, the semilinear Cauchy problem in \eqref{anti deSitter semi} has never been investigated from the viewpoint of blow-up results when \eqref{dominant damping} is satisfied and $m^2\geqslant 0$. As we are going to explain in Subsection \ref{Subsection Main results}, in the case of anti-de Sitter spacetime the growth assumptions on $\Gamma$ depend strongly on the dimension $n$. On the one hand, for low dimensions the situation is quite similar to the corresponding case with exponentially decreasing speed of propagation. On the other hand, for high dimensions the influence of the nonlinear term is dominant and, in particular, the treatment of a threshold case, which could be considered as a critical case in some sense, is more delicate and  requires a more deep analysis of the growth properties of the spatial average of a local solution.

In this last part of the introduction, we recall other known results from the literature on semilinear wave models in de Sitter spacetime and how our result can be framed and understood in relation to these.
Over the last decade several results were established for wave models  by Yagdjian and Yagdjian-Galstian in de Sitter spacetime \cite{YagGal09,Yag09,Yag10,Yag12,Yag13,Yag15,Yag19} with normalized constants $c,H$ (meaning $a_{\mathrm{dS}}(t)=\mathrm{e}^{-t}$) and in anti-de Sitter spacetime (that is, when the speed of propagation is $a_{\mathrm{AdS}}(t)=\mathrm{e}^t$) \cite{YagGal08,YagGal09a,YagGal09e}, respectively. Truly remarkable integral representation formulas for the solutions of the linear Cauchy problem associated with Klein-Gordon equations, both with pure imaginary and real positive mass term, are derived in the case of de Sitter spacetime \cite{YagGal09,Yag09} and anti-de Sitter spacetime \cite{YagGal09a}, respectively. These integral representation formulas have been applied, among other things, to study $L^p-L^q$ estimates, the existence of self-similar solutions, blow-up results with nonlocal nonlinear term as in \eqref{nonlinearity} and to investigate under which assumptions on the coefficients for the mass term and on the space dimension a Huygens' type principle holds. Afterwards, the Cauchy problem associated with the  semilinear wave equation in de Sitter spacetime with power nonlinearity was studied by Nakamura \cite{Na14,Na15} and Ebert-Reissig \cite{ER18} and several global existence results were established not only in classical energy space but also in Sobolev space on $L^2$ basis with different regularities (both below and above the regularity of energy solutions). We point out that in \cite{Na14,Na15} also a nonlinearity of exponential type is considered besides the power nonlinearity. In spite of the above quoted global existence results for small data solutions with a nonnegative power nonlinearity, it seems that there is a lack of understanding concerning the expression for the critical exponent, due to the absence of  a corresponding blow-up counterpart. In this scenario, our results for \eqref{deSitter semi} should emphasize how the presence of the mass term does not allow to prove the blow-up in finite time of any local in time solution when $f(u)=|u|^p$. Indeed, when $\beta=0$ and as $r\to 0^+$ the method that we are going to employ for studying the blow-up is no longer efficient, meaning that the argument that provides the blow-up of the space average fails, with a unique remarkable exception given by the massless case $m^2=0$ (established for the first time in \cite{TsuWa21}, as mentioned above).

\subsection{Main results} \label{Subsection Main results}

Before stating our main results, we introduce the class of solutions to \eqref{deSitter semi} that we will consider throughout this paper.
We emphasize that, even though we will call these solutions weak solutions, we require more regularity than usual distributional solutions. More precisely, we consider the larger class of solutions that can be considered with our approach, and this requires some regularity with respect to the time-variable according to the next definition.

\begin{definition} \label{Def sol} Let $u_0,u_1\in L^1_{\mathrm{loc}}(\mathbb{R}^n)$ such that $\supp u_0, \ \supp u_1 \subset B_R$ for some $R>0$. We say that $$u\in  \mathcal{C}^1\left([0,T), L^{1}_{\mathrm{loc}}(\mathbb{R}^n)\right) \ \mbox{such that} \  f(t,u) \in L^1_{\mathrm{loc}}((0,T)\times \mathbb{R}^n),$$ where the definition of the nonlinear term $f(t,u)$ is given in \eqref{nonlinearity},  is a \emph{weak solution} to \eqref{deSitter semi} on $[0,T)$ if  $u$ fulfills the support condition \eqref{support condition sol} and the integral identity
\begin{align}
& \int_{\mathbb{R}^n} \partial_t u(t,x) \varphi(t,x) \, \mathrm{d}x-\int_{\mathbb{R}^n} u(t,x) \varphi_t(t,x) \, \mathrm{d}x+ b\int_{\mathbb{R}^n}  u(t,x) \varphi(t,x) \, \mathrm{d}x \notag \\
& \qquad + \int_0^t\int_{\mathbb{R}^n} u(s,x) \left(\varphi_{ss}(s,x)-c^2\mathrm{e}^{-2Hs}\Delta \varphi(s,x) - b \,\varphi_s(s,x) +m^2 \varphi(s,x) \right)  \mathrm{d}x \, \mathrm{d}s \notag \\
& \quad  = \varepsilon \int_{\mathbb{R}^n}  u_1(x) \varphi(0,x) \, \mathrm{d}x +\varepsilon \int_{\mathbb{R}^n} u_0(x) \big(b\, \varphi(0,x) -\varphi_t(0,x)\big) \, \mathrm{d}x \notag  \\ & \qquad + \int_0^t \Gamma(s) \left(\int_{\mathbb{R}^n}  |u(s,y)|^p  \, \mathrm{d}y\right)^\beta \int_{\mathbb{R}^n}  |u(s,x)|^p  \varphi(s,x) \, \mathrm{d}x \, \mathrm{d}s  \label{def weak sol int rel 2}
\end{align}
holds for any $t\in (0,T)$ and any test function $\varphi\in\mathcal{C}^\infty_0 ([0,T)\times \mathbb{R}^n)$.
\end{definition}

Considering the following family of time-dependent factors for the nonlinear term in \eqref{nonlinearity}
\begin{align}\label{def Gamma}
\Gamma(t)\doteq \mu \mathrm{e}^{r t} (1+t)^{\kappa}
\end{align}
depending on the real parameters $r,\kappa$, where the multiplicative constant $\mu$ is positive, we are interested in describing how the ranges for $r,\kappa$ affect the blow-up in finite time of local solutions to \eqref{deSitter semi}. In particular, introducing the threshold values
\begin{align}
r_{\mathrm{crit}}(b,m^2,\beta,p) &\doteq \frac{1}{2}\left(b-\sqrt{b^2-4m^2}\right)((\beta+1)p-1), \label{def r crit} \\
\kappa_{\mathrm{crit}}(b,m^2,\beta,p)  &\doteq \begin{cases} -1 & \mbox{if} \ \ b^2>4m^2, \\ -1-(\beta+1)p & \mbox{if} \ \ b^2=4m^2, \end{cases} \label{def kappa crit}
\end{align}
we may distinguish between three different subcases, depending on the range for the parameters $r,\kappa$ in \eqref{def Gamma}:
\begin{itemize}
\item the case with \emph{exponential growth} when $r>r_{\mathrm{crit}}(b,m^2,\beta,p) $ and $\kappa\in\mathbb{R}$;
\item the case with \emph{polynomial growth} when $r=r_{\mathrm{crit}}(b,m^2,\beta,p) $ and $\kappa>\kappa_{\mathrm{crit}}(b,m^2,\beta,p) $;
\item the case with \emph{logarithmic growth} when $r=r_{\mathrm{crit}}(b,m^2,\beta,p) $ and $\kappa=\kappa_{\mathrm{crit}}(b,m^2,\beta,p) $.
\end{itemize} Note that the word ``growth'' in the previous list of subcases does not refer to the growth rate for the function $\Gamma(t)$, rather to the growth of the lower bound for a time-dependent functional related to a local solution $u$, whose evolution in time will be investigated to prove the blow-up in finite time of $u$.

We shall see that a suitable iteration argument for $U$ in Theorem \ref{Thm exponential growth} and for $\mathcal{U}$ in Theorems \ref{Thm polynomial growth} and \ref{Thm logarithmic growth} can be used together with a slicing procedure. For the definition of $U$ and $\mathcal{U}$, respectively, see \eqref{def space average} and \eqref{def mathcal U} below. In the first case we deal with exponential factors in the lower bounds for $U$, while in the threshold case $r=r_{\mathrm{crit}}(b,m^2,\beta,p) $, depending on whether $\kappa>\kappa_{\mathrm{crit}}(b,m^2,\beta,p) $ or $\kappa=\kappa_{\mathrm{crit}}(b,m^2,\beta,p) $, we find lower bounds of polynomial or logarithmic type for $\mathcal{U}$, respectively.

The first result concerns the case with exponential growth.

\begin{theorem}\label{Thm exponential growth} Let $n\geqslant 1$ and $b,m^2\geqslant 0$ such that \eqref{dominant damping} is fulfilled. Let us assume $\beta\geqslant 0, p>1$ and $r>r_{\mathrm{crit}}(b,m^2,\beta,p)$, where $r_{\mathrm{crit}}(b,m^2,\beta,p)$ is defined in \eqref{def r crit}, and consider 
\begin{align}\label{Gamma factor exponential growth}
\Gamma(t)\doteq \mu\,\mathrm{e}^{rt} (1+t)^\kappa
\end{align} for some $\mu>0$ and some $\kappa\in\mathbb{R}$ in \eqref{nonlinearity}. 

 Let us assume that $u_0,u_1\in L^1_{\mathrm{loc}}(\mathbb{R}^n)$ are nonnegative, nontrivial and compactly supported functions with supports contained in $ B_R$ for some $R>0$.

 Let $u\in \mathcal{C}^1\left([0,T), L^{1}_{\mathrm{loc}}(\mathbb{R}^n)\right)$ be a weak solution to the Cauchy problem  \eqref{deSitter semi} according to Definition \ref{Def sol} with lifespan $T=T(\varepsilon)$.
 
 Then, there exists a positive constant $\varepsilon_0=\varepsilon_0(n,c,H,b,m^2,\beta,p,\mu,r,\kappa,u_0,u_1,R)$ such that for any $\varepsilon\in (0,\varepsilon_0]$ the weak solution $u$ blows up in finite time. Furthermore, the following upper bound estimates for the lifespan hold
\begin{align} \label{upper bound for the lifespan exponential growth}
\theta_{b,m^2,p,\beta,r,\kappa} (T(\varepsilon)) \leqslant C \varepsilon^{-\left(\frac{r}{(\beta+1)p-1}-\frac{1}{2}\left(b-\sqrt{b^2-4m^2}\right)\right)^{-1}}, 
\end{align} where the positive constant $C$ is independent of $\varepsilon$ and 
\begin{align}\label{def theta function}
\theta_{b,m^2,p,\beta,r,\kappa}(\tau)\doteq \begin{cases} e^\tau \tau^{\tfrac{\kappa}{r-r_{\mathrm{crit}}(b,m^2,\beta,p)}} & \mbox{if} \ \ b^2>4m^2, \\
e^\tau \tau^{\tfrac{(\beta+1)p-1+\kappa}{r-r_{\mathrm{crit}}(b,m^2,\beta,p)}} & \mbox{if} \ \ b^2=4m^2.\end{cases}
\end{align}
 \end{theorem}
 
\begin{remark} In Theorem \ref{Thm exponential growth} we note that in the case without any additional polynomial growth (or decay) in $\Gamma$, that is for $\kappa=0$, the function in \eqref{def theta function} is simply the exponential function for $b^2>4m^2$. On the other hand,  for $\kappa=0$ and $b^2=4m^2$ we still have a polynomial correction in \eqref{def theta function} (that improves the upper bound estimate).
\end{remark}

The second result concerns the case with polynomial growth. In particular, we have a limit value for the coefficient in the exponential term in \eqref{def Gamma}, while for the polynomial factor we consider the parameter $\kappa$ above the threshold value $\kappa_{\mathrm{crit}}(b,m^2,\beta,p)$.

\begin{theorem}\label{Thm polynomial growth} Let $n\geqslant 1$ and $b,m^2\geqslant 0$ such that \eqref{dominant damping} is fulfilled. Let us assume $\beta\geqslant 0, p>1$ and $r=r_{\mathrm{crit}}(b,m^2,\beta,p)$, $\kappa>\kappa_{\mathrm{crit}}(b,m^2,\beta,p)$, where $r_{\mathrm{crit}}(b,m^2,\beta,p)$ and $\kappa_{\mathrm{crit}}(b,m^2,\beta,p)$ are defined in \eqref{def r crit} and \eqref{def kappa crit}, respectively, and consider 
\begin{align}\label{Gamma factor polynomial growth}
\Gamma(t)\doteq \mu\, \mathrm{e}^{r_{\mathrm{crit}}(b,m^2,\beta,p)t}(1+t)^{\kappa}
\end{align} for some $\mu>0$ in \eqref{nonlinearity}. 

 Let us assume that $u_0,u_1\in L^1_{\mathrm{loc}}(\mathbb{R}^n)$ are nonnegative, nontrivial and compactly supported functions with supports contained in $ B_R$ for some $R>0$.

 Let $u\in \mathcal{C}^1\left([0,T), L^{1}_{\mathrm{loc}}(\mathbb{R}^n)\right)$ be a weak solution to the Cauchy problem  \eqref{deSitter semi} according to Definition \ref{Def sol} with lifespan $T=T(\varepsilon)$.
 
 Then, there exists a positive constant $\varepsilon_0=\varepsilon_0(n,c,H,b,m^2,\beta,p,\mu,\kappa,u_0,u_1,R)$ such that for any $\varepsilon\in (0,\varepsilon_0]$ the weak solution $u$ blows up in finite time. Furthermore, the following upper bound estimates for the lifespan hold
\begin{align} \label{upper bound for the lifespan polynomial growth}
 T(\varepsilon) \leqslant  \begin{cases} C \varepsilon^{-\frac{(\beta+1)p-1}{\kappa+1}} & \mbox{if} \ \, b^2>4m^2, \\
C\varepsilon^{-\left(\frac{\kappa+2}{(\beta+1)p-1}+1\right)^{-1}} & \mbox{if} \ \, b^2= 4m^2,
\end{cases}
\end{align} where the positive constant $C$ is independent of $\varepsilon$.
\end{theorem}

\begin{remark} Let us point out that in the massless case, i.e. for $m^2=0$, we have $r_{\mathrm{crit}}(b,0,\beta,p)=0$. Therefore, in this special case we obtained the blow-up in finite time of local solutions even without requiring additional exponential or polynomial growth for the nonlinear term (setting $\Gamma(t)=1$). In particular, for the local case $\beta=0$, which corresponds to the usual power nonlinearity $|u|^p$, namely, for the semilinear Cauchy problem
\begin{align}\label{deSitter semi power nonlinearity} 
\begin{cases}
\partial_t^2 u- c^2\mathrm{e}^{-2Ht} \Delta u+ b\partial_t u =|u|^p, & x\in \mathbb{R}^n, \ t\in (0,T), \\
u(0,x)= \varepsilon u_0(x), & x\in \mathbb{R}^n, \\
\partial_t u(0,x)= \varepsilon u_1(x), & x\in \mathbb{R}^n,
\end{cases}
\end{align}
 we obtained, with a quite different approach, the same blow-up result recently proved in \cite{TsuWa21}. Moreover, we found the same lifespan estimates as in \cite{TsuWa21}. Indeed, from Theorem \ref{Thm polynomial growth} in this special case the upper bound estimates for the lifespan are given by
 \begin{align*}
 T(\varepsilon) \leqslant  \begin{cases} C \varepsilon^{-(p-1)} & \mbox{if} \ \, b>0, \\
C\varepsilon^{-\frac{p-1}{p+1}} & \mbox{if} \ \, b= 0,
\end{cases}
\end{align*} where the positive constant $C$ is independent of $\varepsilon$.
\end{remark}

The third result concerns the case with logarithmic growth. In this case, we have limit values both for the coefficient of the exponential term and of the polynomial term in \eqref{def Gamma}.

\begin{theorem}\label{Thm logarithmic growth} Let $n\geqslant 1$ and $b,m^2\geqslant 0$ such that \eqref{dominant damping} is fulfilled. Let us assume $\beta\geqslant 0, p>1$ and $r=r_{\mathrm{crit}}(b,m^2,\beta,p)$, $\kappa=\kappa_{\mathrm{crit}}(b,m^2,\beta,p)$, where $r_{\mathrm{crit}}(b,m^2,\beta,p)$ and $\kappa_{\mathrm{crit}}(b,m^2,\beta,p)$ are defined in \eqref{def r crit} and \eqref{def kappa crit}, respectively, and consider 
\begin{align}\label{Gamma factor logarithmic growth}\Gamma(t)\doteq \mu\, \mathrm{e}^{r_{\mathrm{crit}}(b,m^2,\beta,p)t}(1+t)^{\kappa_{\mathrm{crit}}(b,m^2,\beta,p)}
\end{align} for some $\mu>0$ in \eqref{nonlinearity}. 

 Let us assume that $u_0,u_1\in L^1_{\mathrm{loc}}(\mathbb{R}^n)$ are nonnegative, nontrivial and compactly supported functions with supports contained in $ B_R$ for some $R>0$.

 Let $u\in \mathcal{C}^1\left([0,T), L^{1}_{\mathrm{loc}}(\mathbb{R}^n)\right)$ be a weak solution to the Cauchy problem  \eqref{deSitter semi} according to Definition \ref{Def sol} with lifespan $T=T(\varepsilon)$.
 
 Then, there exists a positive constant $\varepsilon_0=\varepsilon_0(n,c,H,b,m^2,\beta,p,\mu,u_0,u_1,R)$ such that for any $\varepsilon\in (0,\varepsilon_0]$ the weak solution $u$ blows up in finite time. Furthermore, the following upper bound estimate for the lifespan holds
\begin{align} \label{upper bound for the lifespan logarithmic growth}
 T(\varepsilon) \leqslant   C \exp\left(K \varepsilon^{-\left((\beta+1)p-1\right)}\right) ,
\end{align} where the positive constants $C,K$ are independent of $\varepsilon$.
\end{theorem}

\begin{remark} \label{Remark improvement Yagjian's result} We emphasize that the results from Theorems \ref{Thm exponential growth}, \ref{Thm polynomial growth} and \ref{Thm logarithmic growth} correspond to ones from Theorem 1.1 in \cite{Yag09}. In particular, in the case $b^2>4m^2$ we improved the limit threshold for the polynomial factor from $\kappa>2$ to $\kappa\geqslant -1$, while we found exactly the same result in the case $b^2=4m^2$. In addition, we established upper bound estimates for the lifespan depending on the precise growth rate of the function $\Gamma(t)$. We underline that in the above mentioned Yagdjian's paper the blow-up of the spatial average of a local solution is proved by means of a comparison argument for certain ordinary differential inequalities, that generalize somehow Kato's lemma (cf. \cite{Kato80}  or \cite{Tak15}). In particular, in that paper, applying the dissipative transformation $w(t,x)=\mathrm{e}^{\frac{b}{2}t}u(t,x)$ and keeping our notations, the equation in \eqref{deSitter semi} is transformed in
\begin{align*}
\partial_t^2 w- c^2\mathrm{e}^{-2Ht} \Delta w-\left(\frac{b^2}{4}-m^2\right) w=\mathrm{e}^{-\frac{b}{2}((\beta+1)p-1)t} f(t,w)
\end{align*} and, then, a modified Kato's lemma is applied to study the blow-up of local in time solutions.
\end{remark}

The second part of the paper will be devoted to the study of blow-up results for local in time solutions to \eqref{anti deSitter semi},
with the time dependent factor $\Gamma$ chosen as follows:
\begin{align}\label{def Gamma anti dS}
\Gamma(t)\doteq \mu \, \mathrm{e}^{\varrho t} (1+t)^{\varsigma}.
\end{align}
 The amplitude function describing the forward light-cone is given by $A_{\mathrm{AdS}}(t)\doteq cH^{-1}(\mathrm{e}^{Ht}-1)$ for this model. 
 
Before stating the main results for \eqref{anti deSitter semi}, also in this case we introduce the class of solutions to \eqref{anti deSitter semi} with which  we will work.

\begin{definition} \label{Def sol anti de Sitter} Let $v_0,v_1\in L^1_{\mathrm{loc}}(\mathbb{R}^n)$ such that $\supp v_0, \ \supp v_1 \subset B_R$ for some $R>0$. We say that $$v\in  \mathcal{C}^1\left([0,T), L^{1}_{\mathrm{loc}}(\mathbb{R}^n)\right) \ \mbox{such that} \  f(t,v) \in L^1_{\mathrm{loc}}((0,T)\times \mathbb{R}^n),$$ where the definition of the nonlinear term $f(t,v)$ is given in \eqref{nonlinearity},  is a \emph{weak solution} to \eqref{anti deSitter semi} on $[0,T)$ if  $v$ fulfills the support condition 
\begin{align} \label{support condition sol anti dS}
\supp v(t,\cdot)\subset B_{R+cH^{-1}(\mathrm{e}^{H t}-1)} \qquad \mbox{for any} \ t\in (0,T),
\end{align} and the integral identity
\begin{align}
& \int_{\mathbb{R}^n} \partial_t v(t,x) \varphi(t,x) \, \mathrm{d}x-\int_{\mathbb{R}^n} v(t,x) \varphi_t(t,x) \, \mathrm{d}x+ b\int_{\mathbb{R}^n}  v(t,x) \varphi(t,x) \, \mathrm{d}x \notag \\
& \qquad + \int_0^t\int_{\mathbb{R}^n} v(s,x) \left(\varphi_{ss}(s,x)-c^2\mathrm{e}^{2Hs}\Delta \varphi(s,x) - b \,\varphi_s(s,x) +m^2 \varphi(s,x) \right)  \mathrm{d}x \, \mathrm{d}s \notag \\
& \quad  = \varepsilon \int_{\mathbb{R}^n}  v_1(x) \varphi(0,x) \, \mathrm{d}x +\varepsilon \int_{\mathbb{R}^n} v_0(x) \big(b\, \varphi(0,x) -\varphi_t(0,x)\big) \, \mathrm{d}x \notag  \\ & \qquad + \int_0^t \Gamma(s) \left(\int_{\mathbb{R}^n}  |v(s,y)|^p  \, \mathrm{d}y\right)^\beta \int_{\mathbb{R}^n}  |v(s,x)|^p  \varphi(s,x) \, \mathrm{d}x \, \mathrm{d}s  \label{def weak sol int rel anti dS}
\end{align}
holds for any $t\in (0,T)$ and any test function $\varphi\in\mathcal{C}^\infty_0 ([0,T)\times \mathbb{R}^n)$.
\end{definition}

Differently from what happens in the case of de Sitter spacetime, when we work in anti-de Sitter spacetime it is possible to derive two different threshold values for the parameter $\varrho$ in \eqref{def Gamma anti dS}. In the next lines we are going to define these two values depending on the range for the space dimension.

We introduce the threshold values 
\begin{align}  \label{def rho crit lin part}
\varrho_{\mathrm{crit}}(n,H,b,m^2,\beta,p)  &\doteq  \frac{1}{2}\left(b-\sqrt{b^2-4m^2}\right)((\beta+1)p-1)+nH(\beta+1)(p-1) 
\end{align} for $n \leqslant\tfrac{\sqrt{b^2-4m^2}}{H} + \tfrac{2}{p}$, and  
\begin{align} \label{def rho crit V0}
\varrho_{\mathrm{crit}}(n,H,b,m^2,\beta,p) &\doteq  \frac{1}{2}(b+nH)((\beta+1)p-1)+nH-(n-1)H(\beta+1)-\frac{H}{p} 
\end{align}   for $n >\tfrac{\sqrt{b^2-4m^2}}{H} + \tfrac{2}{p}$.

The reasons that lead to consider two different values for $\varrho_{\mathrm{crit}}$ depending on whether $n$ is smaller than/equal to or bigger than $\tfrac{\sqrt{b^2-4m^2}}{H} + \tfrac{2}{p}$ and the steps towards to this distinction will be clarified in detail in Subsection \ref{Subsection comparison 1st lb estimates}. Nevertheless, naively and roughly speaking, we can assert that when \eqref{def rho crit lin part} holds the Cauchy data have a stronger influence in the iteration argument than the nonlinear term, while in \eqref{def rho crit V0} the situation is reversed.

The next three theorems are the counterpart in anti-de Sitter spacetime of Theorems \ref{Thm exponential growth}-\ref{Thm logarithmic growth}.

\begin{theorem} \label{Thm anti dS lin lb}
Let $n\geqslant 1$ and $b,m^2\geqslant 0$ such that \eqref{dominant damping} is fulfilled. Let us assume $\beta\geqslant 0$ and  $p>1$ such that 
\begin{align} \label{1st lb for V from lin part is dominant}
\frac{n}{2}-\frac{\sqrt{b^2-4m^2}}{2H}\leqslant \frac{1}{p}
\end{align} and $\varrho>\varrho_{\mathrm{crit}}(n,H,b,m^2,\beta,p)$, where $\varrho_{\mathrm{crit}}(n,H,b,m^2,\beta,p)$ is defined in \eqref{def rho crit lin part}, and consider 
\begin{align*} 
\Gamma(t)\doteq \mu\,\mathrm{e}^{\varrho t} (1+t)^\varsigma
\end{align*} for some $\mu>0$ and some $\varsigma\in\mathbb{R}$ in the term $f(t,v)$ given by \eqref{nonlinearity}. 

 Let us assume that  $v_0,v_1\in L^1_{\mathrm{loc}}(\mathbb{R}^n)$ are nonnegative, nontrivial and compactly supported functions with supports contained in $ B_R$ for some $R>0$.

 Let $v\in \mathcal{C}^1\left([0,T), L^{1}_{\mathrm{loc}}(\mathbb{R}^n)\right)$ be a weak solution to the Cauchy problem  \eqref{anti deSitter semi} according to Definition \ref{Def sol anti de Sitter} with lifespan $T=T(\varepsilon)$.
 
 Then, there exists a positive constant $\varepsilon_0=\varepsilon_0(n,c,H,b,m^2,\beta,p,\mu,\varrho,\varsigma,v_0,v_1,R)$ such that for any $\varepsilon\in (0,\varepsilon_0]$ the weak solution $v$ blows up in finite time. Furthermore, the following upper bound estimates for the lifespan hold
\begin{align*} 
\zeta_{n,H,b,m^2,p,\beta,\varrho,\varsigma} (T(\varepsilon)) \leqslant C  \varepsilon^{-\tfrac{(\beta+1)p-1}{\varrho -\varrho_{\mathrm{crit}}(n,H,b,m^2,\beta,p)}}, 
\end{align*} where the positive constant $C$ is independent of $\varepsilon$ and 
\begin{align*} 
\zeta_{n,H,b,m^2,p,\beta,\varrho,\varsigma}(\tau)\doteq \begin{cases} e^\tau \tau^{\tfrac{\varsigma}{\varrho-\varrho_{\mathrm{crit}}(n,H,b,m^2,\beta,p)}} & \mbox{if} \ \ b^2>4m^2, \\
e^\tau \tau^{\tfrac{(\beta+1)p-1+\varsigma}{\varrho-\varrho_{\mathrm{crit}}(n,H,b,m^2,\beta,p)}} & \mbox{if} \ \ b^2=4m^2.\end{cases}
\end{align*}
\end{theorem}

\begin{theorem} \label{Thm anti dS lin lb poly growth}
Let $n\geqslant 1$ and $b,m^2\geqslant 0$ such that \eqref{dominant damping} is fulfilled. Let us assume $\beta\geqslant 0$ and  $p>1$ satisfying \eqref{1st lb for V from lin part is dominant}  and $\varrho=\varrho_{\mathrm{crit}}(n,H,b,m^2,\beta,p)$ and $\varsigma>\kappa_{\mathrm{crit}}(b,m^2,\beta,p)$, where $\varrho_{\mathrm{crit}}(n,H,b,m^2,\beta,p)$ and $\kappa_{\mathrm{crit}}(b,m^2,\beta,p)$ are defined in \eqref{def rho crit lin part} and in \eqref{def kappa crit}, respectively, and consider 
\begin{align*} 
\Gamma(t)\doteq \mu\,\mathrm{e}^{\varrho_{\mathrm{crit}}(n,H,b,m^2,\beta,p) t} (1+t)^\varsigma
\end{align*} for some $\mu>0$ in the term $f(t,v)$ given by \eqref{nonlinearity}. 

 Let us assume that $v_0,v_1\in L^1_{\mathrm{loc}}(\mathbb{R}^n)$ are nonnegative, nontrivial and compactly supported functions with supports contained in $ B_R$ for some $R>0$.

 Let $v\in \mathcal{C}^1\left([0,T), L^{1}_{\mathrm{loc}}(\mathbb{R}^n)\right)$ be a weak solution to the Cauchy problem  \eqref{anti deSitter semi} according to Definition \ref{Def sol anti de Sitter} with lifespan $T=T(\varepsilon)$.
 
 Then, there exists a positive constant $\varepsilon_0=\varepsilon_0(n,c,H,b,m^2,\beta,p,\mu,\varsigma,v_0,v_1,R)$ such that for any $\varepsilon\in (0,\varepsilon_0]$ the weak solution $v$ blows up in finite time. Furthermore, the following upper bound estimates for the lifespan hold
\begin{align*} 
 T(\varepsilon) \leqslant  \begin{cases} C \varepsilon^{-\frac{(\beta+1)p-1}{\varsigma+1}} & \mbox{if} \ \, b^2>4m^2, \\
C\varepsilon^{-\left(\frac{\varsigma+2}{(\beta+1)p-1}+1\right)^{-1}} & \mbox{if} \ \, b^2= 4m^2,
\end{cases}
\end{align*} where the positive constant $C$ is independent of $\varepsilon$.
\end{theorem}

\begin{theorem} \label{Thm anti dS lin lb log growth}
Let $n\geqslant 1$ and $b,m^2\geqslant 0$ such that \eqref{dominant damping} is fulfilled. Let us assume $\beta\geqslant 0$ and  $p>1$ satisfying \eqref{1st lb for V from lin part is dominant}  and $\varrho=\varrho_{\mathrm{crit}}(n,H,b,m^2,\beta,p)$ and $\varsigma=\kappa_{\mathrm{crit}}(b,m^2,\beta,p)$, where $\varrho_{\mathrm{crit}}(n,H,b,m^2,\beta,p)$ and $\kappa_{\mathrm{crit}}(b,m^2,\beta,p)$ are defined in \eqref{def rho crit lin part} and in \eqref{def kappa crit}, respectively, and consider 
\begin{align*} 
\Gamma(t)\doteq \mu\,\mathrm{e}^{\varrho_{\mathrm{crit}}(n,H,b,m^2,\beta,p) t} (1+t)^{\kappa_{\mathrm{crit}}(b,m^2,\beta,p)}
\end{align*} for some $\mu>0$ in the term $f(t,v)$ given by \eqref{nonlinearity}. 

 Let us assume that $v_0,v_1\in L^1_{\mathrm{loc}}(\mathbb{R}^n)$ are nonnegative, nontrivial and compactly supported functions with supports contained in $ B_R$ for some $R>0$.

 Let $v\in \mathcal{C}^1\left([0,T), L^{1}_{\mathrm{loc}}(\mathbb{R}^n)\right)$ be a weak solution to the Cauchy problem  \eqref{anti deSitter semi} according to Definition \ref{Def sol anti de Sitter} with lifespan $T=T(\varepsilon)$.
 
 Then, there exists a positive constant $\varepsilon_0=\varepsilon_0(n,c,H,b,m^2,\beta,p,\mu,v_0,v_1,R)$ such that for any $\varepsilon\in (0,\varepsilon_0]$ the weak solution $v$ blows up in finite time. Furthermore, the following upper bound estimate for the lifespan holds
 \begin{align*} 
 T(\varepsilon) \leqslant   C \exp\left(K \varepsilon^{-\left((\beta+1)p-1\right)}\right), 
\end{align*} where the positive constants $C,K$ are independent of $\varepsilon$.
\end{theorem}

Let us emphasize that in Theorems \ref{Thm anti dS lin lb}, \ref{Thm anti dS lin lb poly growth} and \ref{Thm anti dS lin lb log growth} we used \eqref{def rho crit lin part} as threshold value for the coefficient $\varrho$ in the exponential term in $\Gamma$, since the condition \eqref{1st lb for V from lin part is dominant}  on the space dimension ensures us that the value for $\varrho_{\mathrm{crit}}$ in \eqref{def rho crit lin part} is smaller than or equal to the one in \eqref{def rho crit V0}. In other words, when \eqref{1st lb for V from lin part is dominant} holds, the wider range for $\varrho$, provided by the condition $\varrho\geqslant \varrho_{\mathrm{crit}}$, is obtained by the definition in \eqref{def rho crit lin part} for $\varrho_{\mathrm{crit}}$.

On the contrary, in the next results we assume that $\varrho_{\mathrm{crit}}$ is given by \eqref{def rho crit V0}, that is, when the opposite inequality of the one in \eqref{1st lb for V from lin part is dominant} holds. 

\begin{theorem} \label{Thm anti dS nlin lb}

Let $n\geqslant 1$ and $b,m^2\geqslant 0$ such that \eqref{dominant damping} is fulfilled. Let us assume $\beta\geqslant 0$ and  $p>1$ such that 
\begin{align}\label{1st lb for V from V0 is dominant}
\frac{n}{2}-\frac{\sqrt{b^2-4m^2}}{2H}>\frac{1}{p}.
\end{align} and $\varrho>\varrho_{\mathrm{crit}}(n,H,b,m^2,\beta,p)$, where $\varrho_{\mathrm{crit}}(n,H,b,m^2,\beta,p)$ is defined in \eqref{def rho crit V0}, and consider 
\begin{align} \label{Gamma factor exponential growth anti dS nlin lb}
\Gamma(t)\doteq \mu\,\mathrm{e}^{\varrho t} (1+t)^\varsigma
\end{align} for some $\mu>0$ and some $\varsigma\in\mathbb{R}$ in the term $f(t,v)$ given by \eqref{nonlinearity}. 

 Let us assume that $v_0,v_1\in L^1_{\mathrm{loc}}(\mathbb{R}^n)$ are nonnegative, nontrivial and compactly supported functions with supports contained in $ B_R$ for some $R>0$.

 Let $v\in \mathcal{C}^1\left([0,T), L^{1}_{\mathrm{loc}}(\mathbb{R}^n)\right)$ be a weak solution to the Cauchy problem  \eqref{anti deSitter semi} according to Definition \ref{Def sol anti de Sitter} with lifespan $T=T(\varepsilon)$.
 
 Then, there exists a positive constant $\varepsilon_0=\varepsilon_0(n,c,H,b,m^2,\beta,p,\mu,\varrho,\varsigma,v_0,v_1,R)$ such that for any $\varepsilon\in (0,\varepsilon_0]$ the weak solution $v$ blows up in finite time. Furthermore, the following upper bound estimates for the lifespan hold
\begin{align} \label{upper bound for the lifespan exponential growth anti dS nlin}
\chi_{n,H,b,m^2,p,\beta,\varrho,\varsigma} (T(\varepsilon)) \leqslant C  \varepsilon^{-\tfrac{(\beta+1)p-1}{\varrho -\varrho_{\mathrm{crit}}(n,H,b,m^2,\beta,p)}}, 
\end{align} where the positive constant $C$ is independent of $\varepsilon$ and 
\begin{align} \label{def chi function anti dS}
\chi_{n,H,b,m^2,p,\beta,\varrho,\varsigma}(\tau)\doteq  e^\tau \tau^{\tfrac{\varsigma}{\varrho-\varrho_{\mathrm{crit}}(n,H,b,m^2,\beta,p)}}. 
\end{align}
\end{theorem}

\begin{remark} Notice that the upper bound estimate for the lifespan \eqref{upper bound for the lifespan exponential growth anti dS nlin} is formally identical to the one in the statement of Theorem \ref{Thm anti dS lin lb} in the case $b^2>4m^2$. Of course, the difference relies in the different definition for the quantity $\varrho_{\mathrm{crit}}$ depending on whether either \eqref{1st lb for V from lin part is dominant} or \eqref{1st lb for V from V0 is dominant} holds.
\end{remark}

\begin{remark} In \cite{PalTak22crit} we provide a blow-up result for \eqref{anti deSitter semi} in the case \eqref{1st lb for V from V0 is dominant} for the critical case $\varrho=\varrho_{\mathrm{crit}}$ by adapting the approach from \cite{TakWak11}.
\end{remark}

\section{Models in de Sitter spacetime}

\subsection{Derivation of the iteration frame}\label{Section iteration frame}

In order to prove Theorems \ref{Thm exponential growth}, we are going to use an iteration argument to show that the space average of a local solution blows up in finite time. Hence, given $u$ local in time solution to \eqref{deSitter semi}, we consider the functional 
\begin{align}\label{def space average}
U(t)\doteq \int_{\mathbb{R}^n} u(t,x) \, \mathrm{d}x \qquad \mbox{for} \ t\in [0,T).
\end{align}
As first step, we are going to determine an iteration frame for the functional $U$. As we explained in the introduction, for us an iteration frame is an integral inequality where $U$ appears both on the left-hand side and on the right-hand side (as a nonlinear term in an integral expression). This iteration frame will allow us to establish a sequence of lower bound estimates of exponential type for $r>r_{\mathrm{crit}}(b,m^2,\beta,p)$, through which we will prove the blow-up in finite time of $U$. 

On the other hand, for the proofs of Theorems \ref{Thm polynomial growth} and \ref{Thm logarithmic growth} rather than with $U(t)$ we will work with the functional $\mathcal{U}(t)$ given by the product of $U(t)$ with a suitable $t$-dependent exponential factor $\mathrm{e}^{\alpha t}$. From the iteration argument for $U$ we will establish immediately the corresponding one for $\mathcal{U}$. By working with $\mathcal{U}$ we will be able to balance the effect of the exponential term in \eqref{def Gamma} in a much more simpler way when $r=r_{\mathrm{crit}}(b,m^2,\beta,p)$. As a result of this balance we may apply a very precise slicing procedure in the different settings of Theorems \ref{Thm polynomial growth} and \ref{Thm logarithmic growth}, depending on whether we work with exponential and/or logarithmic factors and on how many steps are necessary in the slicing procedure.

We point out that in the iteration frame for $U$ (or for $\mathcal{U}$) it is necessary to deal with unbounded exponential multipliers (see also the series of papers \cite{CP19,CP19d,ChenRei20,Chen20,Pal21,Chen20N,PalTak21}, where iteration frames with unbounded exponential multipliers are employed). For this purpose, we apply a \emph{slicing procedure} while deriving the sequence of lower bound estimates for $U$. This procedure is a variation of the first slicing procedure introduced in \cite{AKT00} for the treatment of the critical case for the weakly coupled system of semilinear wave equations in the three dimensional case, where this technique is used to handle factors of logarithmic type. Clearly, the choice of the coefficients characterizing the slicing procedure (see the sequence $\{L_j\}_{j\in\mathbb{N}}$ defined below) is done in order to handle exponential factors in the iteration frame. We will also see how the number of exponential multipliers in the iteration frame will influence the number of steps for the slicing procedure (either a 1 step or a 2 steps procedure).

 Fixed $t\in (0,T)$, we choose a bump function $\varphi\in \mathcal{C}^{\infty}([0,T)\times \mathbb{R}^n)$ that localizes the support of $u$ on the strip $[0,t]\times \mathbb{R}^n$, that is, $\varphi=1$ on $\{(s,x)\in [0,t]\times \mathbb{R}^n:|x|\leqslant R+cH^{-1}(1-\mathrm{e}^{-sH})\}$. Hence, using this $\varphi$ in \eqref{def weak sol int rel 2}, we get
\begin{align*}
& \int_{\mathbb{R}^n}  \partial_t u(t,x) \mathrm{d}x+b \int_{\mathbb{R}^n}  u(t,x) \mathrm{d}x  +m^2\int_0^t  \int_{\mathbb{R}^n} u(s,x) \mathrm{d}x \, \mathrm{d}s  \\ & \qquad = \varepsilon \int_{\mathbb{R}^n} (u_1(x)+bu_0(x)) \, \mathrm{d}x +\int_0^t \Gamma(s)\left( \int_{\mathbb{R}^n}  |u(s,x)|^p \mathrm{d}x\right)^{\beta+1}  \mathrm{d}s,
\end{align*} that can be rewritten as
\begin{align*}
U'(t)+bU(t)+m^2 \int^t_0U(s) \mathrm{d}s = \varepsilon \int_{\mathbb{R}^n} (u_1(x)+bu_0(x)) \, \mathrm{d}x +\int_0^t \Gamma(s)\left( \int_{\mathbb{R}^n}  |u(s,x)|^p \mathrm{d}x\right)^{\beta+1}  \mathrm{d}s.
\end{align*} From the previous relation we see that $U$ is twice continuously differentiable and that
\begin{align}\label{ODE for F}
U''(t)+bU'(t)+m^2U(t)= \Gamma(t)\left( \int_{\mathbb{R}^n}  |u(t,x)|^p \mathrm{d}x\right)^{\beta+1}.
\end{align} 
Thanks to the assumption \eqref{dominant damping}, we may factorize the differential operator on the left-hand side of \eqref{ODE for F} as follows:
\begin{align}\label{decomposition 2nd order operator}
\mathrm{e}^{-\alpha_1 t} \frac{\mathrm{d}}{\mathrm{d} t} \left( \mathrm{e}^{(\alpha_1-\alpha_2) t} \frac{\mathrm{d}}{\mathrm{d} t} \left( \mathrm{e}^{\alpha_2 t}U(t)\right) \right) =  U''(t)+(\alpha_1+\alpha_2)U'(t)+\alpha_1 \alpha_2 U(t),
\end{align}
where the pair of real parameters $(\alpha_1,\alpha_2)$ satisfies
\begin{align*}
\alpha_1+\alpha_2=b, \quad \alpha_1\alpha_2= m^2.
\end{align*} 
 Clearly, the previous conditions for $\alpha_1$ and $\alpha_2$ are symmetric and they are satisfied if $\alpha_{1/2}$ are the roots of the quadratic equation 
 \begin{align}\label{quadratic eq gamma}
 \alpha^2-b\alpha+m^2=0.
 \end{align} Note that in the balanced case $b^2=4m^2$ the previous equation has a double root and $\alpha_1=\alpha_2= \frac{b}{2}$ and that in the dominant mass case the roots of \eqref{quadratic eq gamma} are complex conjugate, so oscillations appear.
Therefore, we may rewrite \eqref{ODE for F} as follows:
\begin{align}\label{ODE for F v2}
\mathrm{e}^{-\alpha_1 t} \frac{\mathrm{d}}{\mathrm{d} t} \left( \mathrm{e}^{(\alpha_1-\alpha_2) t} \frac{\mathrm{d}}{\mathrm{d} t} \left( \mathrm{e}^{\alpha_2 t}U(t)\right) \right) = \Gamma(t)\left( \int_{\mathbb{R}^n}  |u(t,x)|^p \mathrm{d}x\right)^{\beta+1}.
\end{align} 

Next, we can use \eqref{ODE for F v2} to derive the iteration frame for $U$ by assuming nonnegative $u_0$ and $u_1$. Let us begin with the case $b^2>4m^2$ (when $\alpha_1\neq \alpha_2$). Multiplying \eqref{ODE for F v2} by $\mathrm{e}^{\alpha_1 t}$ and integrating over $[0,t]$, we find
\begin{align*}
\int_0^t \mathrm{e}^{\alpha_1 \tau} \,\Gamma(\tau)\left( \int_{\mathbb{R}^n}  |u(\tau,x)|^p \mathrm{d}x\right)^{\beta+1} \mathrm{d}\tau = \mathrm{e}^{(\alpha_1-\alpha_2) t} \frac{\mathrm{d}}{\mathrm{d} t} \left( \mathrm{e}^{\alpha_2 t}U(t)\right) -(U'(0)+\alpha_2 U(0)).
\end{align*} Analogously, from this last relation we obtain
\begin{align*}
& \int_0^t \mathrm{e}^{(\alpha_2-\alpha_1) s} \int_0^s \mathrm{e}^{\alpha_1 \tau} \,\Gamma(\tau)\left( \int_{\mathbb{R}^n}  |u(\tau,x)|^p \mathrm{d}x\right)^{\beta+1} \mathrm{d}\tau \, \mathrm{d}s \\ & \qquad =  \mathrm{e}^{\alpha_2 t}U(t) - U(0) +\frac{\mathrm{e}^{(\alpha_2-\alpha_1) t} -1}{\alpha_1-\alpha_2}  (U'(0)+\alpha_2 U(0)),
\end{align*} which implies in turn
\begin{align}
U(t) & = \frac{\alpha_2\, \mathrm{e}^{-\alpha_1 t}-\alpha_1\,\mathrm{e}^{-\alpha_2 t}}{\alpha_2-\alpha_1}U(0) + \frac{\mathrm{e}^{-\alpha_1 t}-\mathrm{e}^{-\alpha_2 t}}{\alpha_2-\alpha_1}U'(0)  \notag \\ & \qquad  + \mathrm{e}^{-\alpha_2 t} \int_0^t \mathrm{e}^{(\alpha_2-\alpha_1) s} \int_0^s \mathrm{e}^{\alpha_1 \tau}\, \Gamma(\tau)\left( \int_{\mathbb{R}^n}  |u(\tau,x)|^p \mathrm{d}x\right)^{\beta+1} \mathrm{d}\tau \, \mathrm{d}s \notag \\
& = \varepsilon \, \frac{\alpha_2\, \mathrm{e}^{-\alpha_1 t}-\alpha_1\,\mathrm{e}^{-\alpha_2 t}}{\alpha_2-\alpha_1} \int_{\mathbb{R}^n} u_0(x) \, \mathrm{d} x +  \varepsilon \,\frac{\mathrm{e}^{-\alpha_1 t}-\mathrm{e}^{-\alpha_2 t}}{\alpha_2-\alpha_1}\int_{\mathbb{R}^n} u_1(x) \, \mathrm{d} x  \notag \\ & \qquad +  \mathrm{e}^{-\alpha_2 t} \int_0^t \mathrm{e}^{(\alpha_2-\alpha_1) s} \int_0^s \mathrm{e}^{\alpha_1 \tau} \, \Gamma(\tau)\left( \int_{\mathbb{R}^n}  |u(\tau,x)|^p \mathrm{d}x\right)^{\beta+1} \mathrm{d}\tau \, \mathrm{d}s. \label{representation F dom damp case}
\end{align} In the limit case $b^2=4m^2$, we can proceed similarly obtaining
\begin{align}
U(t) &= \mathrm{e}^{-\frac{b}{2} t}U(0) +  \left(U'(0)+\tfrac{b}{2} U(0)\right) t \, \mathrm{e}^{-\frac{b}{2} t}  + \mathrm{e}^{-\frac{b}{2} t} \int_0^t  \int_0^s \mathrm{e}^{\frac{b}{2} \tau} \, \Gamma(\tau)\left( \int_{\mathbb{R}^n}  |u(\tau,x)|^p \mathrm{d}x\right)^{\beta+1} \mathrm{d}\tau \, \mathrm{d}s  \notag \\
& = \varepsilon  \left(1+\tfrac{b}{2}t\right) \mathrm{e}^{-\frac{b}{2} t}\int_{\mathbb{R}^n} u_0(x) \, \mathrm{d} x +  \varepsilon \, t \, \mathrm{e}^{-\frac{b}{2} t} \int_{\mathbb{R}^n} u_1(x) \, \mathrm{d} x  \notag  \\ & \qquad + \mathrm{e}^{-\frac{b}{2}t} \int_0^t \int_0^s \mathrm{e}^{\frac{b}{2} \tau}  \, \Gamma(\tau)\left( \int_{\mathbb{R}^n}  |u(\tau,x)|^p \mathrm{d}x\right)^{\beta+1}\mathrm{d}\tau \, \mathrm{d}s. \label{representation F balanced case}
\end{align}

Consequently, requiring that $u_0$ and $u_1$ are nonnegative functions, then, from the previous identities we obtain immediately that $U$ is a nonnegative functional. Next we determine the iteration frame. Since $\supp u(t,\cdot) \subset B_{R+A_{\mathrm{dS}}(t)}$ for any $t\in (0,T)$, by using H\"older's inequality we have
\begin{align*}
0\leqslant U(t) & \leqslant \left(\int_{\mathbb{R}^n}|u(t,x)|^p \, \mathrm{d}x\right)^{\frac{1}{p}}  \left(\meas\left(B_{R+A_{\mathrm{dS}}(t)}\right)\right)^{\frac{1}{p'}} \\ &\lesssim (R+A_{\mathrm{dS}}(t))^{\frac{n}{p'}} \left(\int_{\mathbb{R}^n}|u(t,x)|^p \, \mathrm{d}x\right)^{\frac{1}{p}} \lesssim \left(\int_{\mathbb{R}^n}|u(t,x)|^p \, \mathrm{d}x\right)^{\frac{1}{p}}, 
\end{align*} and, hence,
\begin{align*}
\int_{\mathbb{R}^n}|u(t,x)|^p \, \mathrm{d}x  \gtrsim  (U(t))^p.
\end{align*} 
Notice that in the previous step, we took advantage of the fact that the light-cone is contained in an infinite half cylinder.

Thus, from \eqref{representation F dom damp case} and \eqref{representation F balanced case} we get the iteration frame
\begin{equation} \label{Iteration frame}
 U(t) \geqslant C  \mathrm{e}^{-\alpha_2 t} \int_0^t \mathrm{e}^{(\alpha_2-\alpha_1) s} \int_0^s \mathrm{e}^{\alpha_1 \tau}  \, \Gamma(\tau) (U(\tau))^{(\beta+1)p} \, \mathrm{d}\tau \, \mathrm{d}s, 
\end{equation}
where $C=C(n,c,H,p,R)>0$ is a suitable constant.

Clearly, in order to be able to apply the previous iteration frames to get a sequence of lower bound estimates for $U$, we need to determine a first lower bound for $U$. From \eqref{representation F dom damp case} and \eqref{representation F balanced case}, since the Cauchy data are taken nonnegative and nontrivial, we have immediately the lower bound estimates
\begin{align} \label{1st lb U}
U(t)\geqslant \begin{cases} K_0 \varepsilon \, \mathrm{e}^{-\left(\frac{b}{2}-\frac{1}{2}\sqrt{b^2-4m^2}\right)t} & \mbox{if} \ b^2>4m^2, \\
K_0 \varepsilon \, (1+t) \, \mathrm{e}^{-\frac{b}{2} t} & \mbox{if} \ b^2=4m^2, 
\end{cases} 
\end{align} for any $t\in (0,T)$, where $K_0=K_0(b,m^2,u_0,u_1)$ is a suitable positive and independent of $\varepsilon$ constant.

We emphasize that \eqref{Iteration frame} is the iteration frame that will be used in the proof of Theorem \ref{Thm exponential growth}, while for Theorems \ref{Thm polynomial growth} and \ref{Thm logarithmic growth} the choice of the time-dependent functional and the corresponding iteration frame will follow directly from \eqref{Iteration frame}.

In the next three subsections, we will prove these theorems. In each case the growth condition assumed on $\Gamma$ has a crucial role in determining the key factors in the iteration frame and, consequently, the main features of the associated slicing procedure.

\subsection{Case with exponential growth: proof of Theorem \ref{Thm exponential growth}} \label{Subsection exponential growth}

In this subsection, we prove Theorem \ref{Thm exponential growth}. As anticipated, the time-dependent functional that we consider to prove the blow-up result is the space average $U$ defined in \eqref{def space average}.

Since the time-dependent factor $\Gamma$ in \eqref{nonlinearity} is given by \eqref{Gamma factor exponential growth} with $r>r_{\mathrm{crit}(b,m^2,\beta,p)}$, the exponential growth of $\Gamma$ is dominant over the first lower bound for $U$ in \eqref{1st lb U} (which decays exponentially). Therefore, when deriving the sequence of lower bound estimates for $U$ through \eqref{Iteration frame}, we need to handle exponentially increasing factors both in the $\tau$-integral and in the $s$-integral. Hence, we apply a 2 steps slicing procedure and the coefficients characterizing the shrinking of the domains of integration on the right-hand side of \eqref{Iteration frame} are chosen in order to allow the handling of the unbounded exponential multipliers $\mathrm{e}^{(\alpha_1+r)\tau}$ and $\mathrm{e}^{(\alpha_2+r)s}$ in the first and in the second integral, respectively. 

We may define now the parameters $\{L_j\}_{j\in\mathbb{N}}$ that characterize the slicing procedure:
\begin{align}\label{def Lj exponential growth}
L_j\doteq \prod_{k=0}^j \ell_k \qquad \mbox{for any} \ j\in\mathbb{N},
\end{align} where the coefficients $\{\ell_k\}_{k\in\mathbb{N}}$ are given by
\begin{align*}
\ell_0 & \doteq \max\left\{(r+\alpha_1)^{-1},(r+\alpha_2)^{-1}\right\}, \qquad
\ell_k  \doteq 1+\left((\beta+1)p\right)^{-k/2} \quad \mbox{for any} \ k\geqslant 1.
\end{align*}
Notice that $\ell_0$ is well defined thanks to the condition on $r$. Moreover, since $\ell_k>1$ for any $k\geqslant 1$, the sequence $\{L_j\}_{j\in\mathbb{N}}$ is strictly increasing. Finally, due to the choice of  $\{\ell_k\}_{k\geqslant 1}$, we have that the series $\sum_{k=1}^{\infty}\ln \ell_k$ is convergent, and this is equivalent to prove the convergence of the following infinite product $$L\doteq \prod_{j=0}^\infty L_j \in\mathbb{R}_+ .$$
Our first goal is to prove the following sequence of lower bound estimates for $U$:
\begin{align}\label{lb estimate U step j exponential growth}
U(t)\geqslant C_j \mathrm{e}^{a_j t}(t-L_{2j})^{b_j} (1+t)^{-\beta_j} \qquad \mbox{for} \ t\geqslant L_{2j} \ \mbox{and for any} \ j\in\mathbb{N},
\end{align} where $\{C_j\}_{j\in\mathbb{N}}$, $\{a_j\}_{j\in\mathbb{N}\setminus\{0\}}$, $\{b_j\}_{j\in\mathbb{N}}$, $\{\beta_j\}_{j\in\mathbb{N}}$ are sequences of nonnegative real numbers to be determined iteratively.

For $j=0$ \eqref{lb estimate U step j exponential growth} is given by \eqref{1st lb U} provided that $C_0\doteq K_0\varepsilon$, $a_0\doteq -\frac{b}{2}+\frac{1}{2}\sqrt{b^2-4m^2}$, $\beta_0=0$, and, finally, $b_0\doteq0$ if $b^2>4m^2$ and  $b_0\doteq 1$ if $b^2=4m^2$. We underline that $a_0$ is the only term in the sequence $\{a_j\}_j$ that is not positive.

Denoting by $\kappa_+$ and $\kappa_-$ the positive and the negative part of $\kappa$ (i.e., $\kappa_+\doteq\max\{\kappa,0\}$ and $\kappa_-\doteq - \min\{	\kappa,0\}$), from \eqref{Iteration frame} we get 
\begin{equation} \label{Iteration frame exponential growth}
 U(t) \geqslant \mu C (1+t)^{-\kappa_-} \mathrm{e}^{-\alpha_2 t} \int_0^t \mathrm{e}^{(\alpha_2-\alpha_1) s} \int_0^s \mathrm{e}^{(\alpha_1+r) \tau} \tau^{\kappa_+} (U(\tau))^{(\beta+1)p} \, \mathrm{d}\tau \, \mathrm{d}s.
\end{equation}
We want to prove \eqref{lb estimate U step j exponential growth} by induction with respect to $j$. We have already remarked the validity of the base case. Next we prove the induction step. Assuming that \eqref{lb estimate U step j exponential growth} is satisfied for some $j\geqslant 0$ we prove it for $j+1$. Plugging the lower bound estimate \eqref{lb estimate U step j exponential growth} in \eqref{Iteration frame exponential growth}, for $t\geqslant L_{2j}$ we have
\begin{align*}
U(t) & \geqslant \mu C (1+t)^{-\kappa_-} \mathrm{e}^{-\alpha_2 t} \int_{L_{2j}}^t \mathrm{e}^{(\alpha_2-\alpha_1) s} \int_{L_{2j}}^s \mathrm{e}^{(\alpha_1+r) \tau} \tau^{\kappa_+} (U(\tau))^{(\beta+1)p} \, \mathrm{d}\tau \, \mathrm{d}s \\
& \geqslant \mu C  C_j^{q} (1+t)^{-(\kappa_- +q\beta_j) } \mathrm{e}^{-\alpha_2 t} \int_{L_{2j}}^t \mathrm{e}^{(\alpha_2-\alpha_1) s} \int_{L_{2j}}^s \mathrm{e}^{(\alpha_1+r+q a_j) \tau} (\tau-L_{2j})^{\kappa_++qb_j} \, \mathrm{d}\tau \, \mathrm{d}s,
\end{align*} where from now on, for the sake of brevity, we denote $q\doteq (\beta+1)p$. For $t\geqslant L_{2j+1}$ we can shrink the domain of integration in the previous inequality as follows:
\begin{align}
U(t) &  \geqslant \mu C  C_j^{q} (1+t)^{-(\kappa_- +q\beta_j) } \mathrm{e}^{-\alpha_2 t} \int_{L_{2j+1}}^t \mathrm{e}^{(\alpha_2-\alpha_1) s} \int_{\tfrac{L_{2j}s}{L_{2j+1}}}^s \mathrm{e}^{(\alpha_1+r+q a_j) \tau} (\tau-L_{2j})^{\kappa_++qb_j} \, \mathrm{d}\tau \, \mathrm{d}s \notag \\
 &  \geqslant \frac{\mu C  C_j^{q}}{\ell_{2j+1}^{\kappa_+ +q b_j }} (1+t)^{-(\kappa_- +q\beta_j) } \mathrm{e}^{-\alpha_2 t} \!\int_{L_{2j+1}}^t \! \mathrm{e}^{(\alpha_2-\alpha_1) s} (s-L_{2j+1})^{\kappa_++qb_j} \! \int_{\tfrac{s}{\ell_{2j+1}}}^s \!\mathrm{e}^{(\alpha_1+r+q a_j) \tau}  \mathrm{d}\tau \, \mathrm{d}s, \label{intermediate lb U 1st step exponential case}
\end{align} where in the second step we used the monotonicity of the factor $(\tau-L_{2j})^{\kappa_++qb_j}$. Let us show now how we can estimate from below the $\tau$-integral. By a direct computation we have
\begin{align*}
\int_{\tfrac{s}{\ell_{2j+1}}}^s \mathrm{e}^{(\alpha_1+r+q a_j) \tau} \, \mathrm{d}\tau & = (\alpha_1+r+q a_j)^{-1}  \mathrm{e}^{(\alpha_1+r+q a_j) s} \left(1- \mathrm{e}^{-(\alpha_1+r+q a_j)(1-1/\ell_{2j+1}) s}\right) \\
 & \geqslant (\alpha_1+r+q a_j)^{-1}  \mathrm{e}^{(\alpha_1+r+q a_j) s} \left(1- \mathrm{e}^{-(\alpha_1+r+q a_j)(\ell_{2j+1}-1) L_{2j}}\right) \\
 & \geqslant (\alpha_1+r+q a_j)^{-1}  \mathrm{e}^{(\alpha_1+r+q a_j) s} \left(1- \mathrm{e}^{-(\alpha_1+r)(\ell_{2j+1}-1) L_{0}}\right)  \\
 & \geqslant (\alpha_1+r+q a_j)^{-1}  \mathrm{e}^{(\alpha_1+r+q a_j) s} \left(1- \mathrm{e}^{-(\ell_{2j+1}-1)}\right)
\end{align*} for $s\geqslant L_{2j+1}$, where in the previous chain of inequalities we used the following properties $\alpha_1+r>0$, $a_j\geqslant 0$, $\ell_{2j+1}>1$, $L_{2j}\uparrow$ and $L_0=\ell_0\geqslant (\alpha_1+r)^{-1}$. Then, using the inequality $1-\mathrm{e}^{-y}\geqslant y-\frac{y^2}{2}$ for any $y\geqslant 0$ we have
\begin{align}
1- \mathrm{e}^{-(\ell_{2j+1}-1)} &\geqslant (\ell_{2j+1}-1) \left(1-\tfrac{1}{2}(\ell_{2j+1}-1)\right) = q^{-(2j+1)}(q^{j+1/2}-\tfrac{1}{2}) \notag \\
&\geqslant q^{-(2j+1)}(q-\tfrac{1}{2}). \label{remainder exponential term}
\end{align} Combining these two last inequalities, from \eqref{intermediate lb U 1st step exponential case} we obtain
\begin{align*}
U(t)  &  \geqslant \frac{\mu C (q-\tfrac 12) C_j^{q}}{\ell_{2j+1}^{\kappa_+ +q b_j } (\alpha_1+r+q a_j) q^{2j+1}} (1+t)^{-(\kappa_- +q\beta_j) } \mathrm{e}^{-\alpha_2 t} \! \int_{L_{2j+1}}^t \!  \mathrm{e}^{(\alpha_2+r+q a_j) s} (s-L_{2j+1})^{\kappa_++qb_j}  \mathrm{d}s.
\end{align*}
Until now we applied a first step in the slicing procedure to deal with the $\tau$-integral. Repeating analogous computations after shrinking the domain of integration to $[t/\ell_{2j+2},t]$ in the $s$-integral for $t\geqslant L_{2j+2}$, we arrive at the lower bound estimate
\begin{align*}
U(t)  &  \geqslant \frac{\mu C (q-\tfrac 12)^2 C_j^{q}(\ell_{2j+1}\ell_{2j+2})^{-(\kappa_+ +q b_j) }}{ (\alpha_1+r+q a_j) (\alpha_2+r+q a_j) q^{4j+3}} \,  \mathrm{e}^{(r+q a_j) t} (t-L_{2j+2})^{\kappa_++qb_j} (1+t)^{-(\kappa_- +q\beta_j) }, 
\end{align*} which is exactly \eqref{lb estimate U step j exponential growth} for $j+1$ provided that
\begin{align}
C_{j+1} & \doteq \frac{\mu C (q-\tfrac 12)^2 C_j^{q}(\ell_{2j+1}\ell_{2j+2})^{-(\kappa_+ +q b_j) }}{ (\alpha_1+r+q a_j) (\alpha_2+r+q a_j) q^{4j+3}}, \label{recursive relation for Cj} \\
a_{j+1} & \doteq r+q a_j, \quad  b_{j+1}  \doteq \kappa_+ +q b_j, \quad \beta_{j+1}  \doteq \kappa_- +q \beta_j.  \label{recursive relation for aj, bj, beta j} 
\end{align}
By applying recursively the previous relations among two consecutive terms from the sequences $\{a_j\}_{j\in\mathbb{N}}$, $\{b_j\}_{j\in\mathbb{N}}$, $\{\beta_j\}_{j\in\mathbb{N}}$ we obtain the explicit representation
\begin{align}
a_j= r\sum_{k=0}^{j-1} q^k+q^j a_0 = \frac{q^j-1}{q-1} r+q^j a_0 = \left(\frac{r}{q-1}+a_0\right)q^j - \frac{r}{q-1}, \label{representation aj}
\end{align} and, in an analogous way,
\begin{align}
b_j & = \left(\frac{\kappa_+}{q-1}+b_0\right)q^j - \frac{\kappa_+}{q-1}, \label{representation bj} \\
\beta_j & = \frac{\kappa_-}{q-1} q^j - \frac{\kappa_-}{q-1}, \label{representation beta j}
\end{align} where in the last relation we used $\beta_0=0$.
The next step is to determine a lower bound for the constant $C_j$ that we can handle more easily. We remark that, since $r>r_{\mathrm{crit}}(b,m^2,\beta,p)$ the quantity $r/(q-1)+a_0$ is strictly positive. Therefore,
\begin{align}
\alpha_{1/2}+r+q a_j &= \alpha_{1/2}+ a_{j+1} < (\tfrac{r}{q-1}+a_0)q^j+ \alpha_{1/2}+a_0 \notag \\
 & \leqslant  (\tfrac{r}{q-1}+a_0)q^j +\sqrt{b^2-4m^2} \notag \\ 
 & \leqslant M_0 q^j \label{upper bound alpha+r+qa j}
\end{align} for any $j\in\mathbb{N}$, where $M_0=M_0(b,m^2,r,\beta,p)$ is a suitable positive quantity which is independent of $j$.
Furthermore, we remark that 
\begin{align*}
\lim_{j\to \infty} (\ell_{2j+1}\ell_{2j+2})^{b_{j+1} } & = \lim_{j\to \infty} \exp\left(b_{j+1}\left[\ln\ell_{2j+1}+\ln\ell_{2j+2}\right]\right) \\ & = \lim_{j\to \infty} \exp\left(\left(\tfrac{\kappa_+}{q-1}+b_0\right)q^{j+1}\left[\ln \left(1+q^{-(j+1/2)}\right)+\ln\left(1+q^{-(j+1)}\right)\right]\right) \\ & = \exp\left(\left(\tfrac{\kappa_+}{q-1}+b_0\right)(1+\sqrt{q})\right),
\end{align*} consequently, there exists a uniform (i.e. independent of $j$) constant $M_1=M_1(b,m^2,\kappa,\beta,p)>0$ such that $ (\ell_{2j+1}\ell_{2j+2})^{b_{j+1} }\leqslant M_1$ for any $j\in \mathbb{N}$. Combining \eqref{recursive relation for Cj}, \eqref{recursive relation for aj, bj, beta j}, \eqref{upper bound alpha+r+qa j} and the previous uniform upper bound, we see that
\begin{align*}
C_{j+1} & = \frac{\mu C (q-\tfrac 12)^2  C_j^{q}}{ (\ell_{2j+1}\ell_{2j+2})^{b_{j+1}} (\alpha_1+ a_{j+1}) (\alpha_2+a_{j+1}) q^{4j+3}} \geqslant \underbrace{\frac{\mu C (q-\tfrac 12)^2 q^3}{ M_1 M_0^2 }}_{D\doteq } q^{- 6(j+1)}C_j^{q}.
\end{align*} We can now use the inequality $C_j\geqslant D q^{-6j} C_{j-1}^q$ to derive a more convenient lower bound for $C_j$ for sufficiently large indexes. Applying the logarithmic function to both sides of the previous inequality and using iteratively the resulting inequality, we find 
\begin{align*}
\ln C_j &  \geqslant q\ln C_{j-1}-6j \ln q+\ln D \geqslant q^2\ln C_{j-1}-6(j+(j-1)q) \ln q+ (1+q)\ln D \\
 &  \geqslant \ldots \geqslant q^j \ln C_0-6 \left(\,\sum_{k=0}^{j-1}(j-k)q^k\right)\ln q+\left(\,\sum_{k=0}^{j-1}q^k\right)\ln D.
\end{align*} Using the following identity 
\begin{align} \label{summation identity}
\sum_{k=0}^{j-1}(j-k)q^k= \frac{1}{q-1}\left(\frac{q^{j+1}-q}{q-1}-j\right),
\end{align} we have
\begin{align*}
\ln C_j &  \geqslant q^j\left( \ln C_0-\frac{6q \ln q}{(q-1)^2}+\frac{\ln D}{q-1}\right)+\frac{6q \ln q}{(q-1)^2}
+\frac{6\ln q}{q-1} j-\frac{\ln D}{q-1}.
\end{align*} Let $j_0=j_0(n,c,H,b,m^2,\mu,r,\kappa,\beta,p,R)\in\mathbb{N}$ be the smallest integer such that $j_0\geqslant \frac{\ln D}{6\ln q}-\frac{q}{q-1}$. Then, for any $j\geqslant j_0$ it results
\begin{align}\label{exponential lb for Cj}
\ln C_j &  \geqslant q^j\left( \ln (K_0 \varepsilon)-\frac{6q \ln q}{(q-1)^2}+\frac{\ln D}{q-1}\right)= q^j \ln (\tilde{D}\varepsilon),
\end{align} where $\tilde{D}\doteq K_0 q^{-6q/(q-1)^2}D^{1/(q-1)}$.
Hence, recalling that $L_{2j}\uparrow L$, if we combine \eqref{lb estimate U step j exponential growth}, \eqref{representation aj}, \eqref{representation bj}, \eqref{representation beta j} and \eqref{exponential lb for Cj} for $t\geqslant L$ and for any $j\geqslant j_0$ it holds
\begin{align*}
U(t) & \geqslant \exp\left(q^j \left(\ln (\tilde{D}\varepsilon)+\left(\tfrac{r}{q-1}+a_0\right) t+\left(\tfrac{\kappa_+}{q-1}+b_0\right)\ln (t-L)-\tfrac{\kappa_-}{q-1}\ln (1+t) \right)\right) \\ & \qquad \times \exp\left(-\tfrac{r}{q-1} t\right)(t-L)^{-\frac{\kappa_+}{q-1}} (1+t)^{-\frac{\kappa_-}{q-1}}. 
\end{align*} Next, using the trivial inequalities $\ln(t-L)\geqslant \ln t -\ln2$ and $-\ln(1+t)\geqslant -\ln t-\ln2$ for $t\geqslant \max\{2L,1\}$ and the identity $\kappa=\kappa_+-\kappa_-$, from the previous estimate we obtain
\begin{align*}
U(t) & \geqslant \exp\left(q^j \ln \left(\hat{D}\varepsilon \left(t^{\tfrac{\kappa+(q-1)b_0}{r+(q-1)a_0}} \mathrm{e}^t\right)^{\frac{r}{q-1}+a_0}\right)\right)  \exp\left(-\tfrac{r}{q-1} t\right)(t-L)^{-\frac{\kappa_+}{q-1}} (1+t)^{-\frac{\kappa_-}{q-1}},
\end{align*} where $\hat{D}\doteq 2^{-((\kappa_++\kappa_-)/(q-1) +b_0)}\tilde{D}$. By using the function defined in \eqref{def theta function}, we may rewrite 
\begin{align}\label{final lb U exponential growth}
U(t) & \geqslant \exp\left(q^j \ln \left(\hat{D}\varepsilon \left(\theta_{b,m^2,p,\beta,r,\kappa}(t)\right)^{\frac{r}{q-1}+a_0}\right)\right)  \exp\left(-\tfrac{r}{q-1} t\right)(t-L)^{-\frac{\kappa_+}{q-1}} (1+t)^{-\frac{\kappa_-}{q-1}}
\end{align} for $t\geqslant \max\{2L,1\}$ and for $j\geqslant j_0$.

 From \eqref{def theta function} we see that $\theta_{b,m^2,p,\beta,r,\kappa}$ is strictly increasing (and hence invertible) for $t\geqslant \tilde{T}$, where $\tilde{T}=\tilde{T}(b,m^2,p,\beta,r,\kappa)$ is a suitable nonnegative quantity. Note that for $\kappa\geqslant 0$ if $b^2>4m^2$ and $\kappa\geqslant -(\beta+1)p+1$ if $b^2>4m^2$ we can simply take $\tilde{T}(b,m^2,p,\beta,r,\kappa)=0$. With a slight abuse of notation, in what follows we denote by $\theta_{b,m^2,p,\beta,r,\kappa}^{-1}$ the inverse function of the restriction $\theta_{b,m^2,p,\beta,r,\kappa}\big|_{[\tilde{T},\infty)}$.
 
 We remark that the logarithmic factor multiplying $q^j$ in \eqref{final lb U exponential growth} is strictly positive if and only if $\hat{D}\varepsilon \left(\theta_{b,m^2,p,\beta,r,\kappa}(t)\right)^{\frac{r}{q-1}+a_0}>1$. For $t\geqslant \tilde{T}$, this is equivalent to require
 \begin{align*}
 t> \theta_{b,m^2,p,\beta,r,\kappa}^{-1}\left((\hat{D}\varepsilon)^{-\left(\frac{r}{q-1}+a_0\right)^{-1}}\right).
 \end{align*} 
 
 Since $\lim_{s\to \infty} \theta_{b,m^2,p,\beta,r,\kappa}^{-1}(s)=\infty$, we may fix $\varepsilon_0=\varepsilon_0(n,c,H,b,m^2,q,\mu,r,\kappa,u_0,u_1,R)>~0$ sufficiently small so that $$\theta_{b,m^2,p,\beta,r,\kappa}^{-1}\left((\hat{D}\varepsilon_0)^{-\left(\frac{r}{q-1}+a_0\right)^{-1}}\right)\geqslant \max\{2L,1,\tilde{T}\}.$$ Thus, for any $\varepsilon\in(0,\varepsilon_0]$ and any $t>\theta_{b,m^2,p,\beta,r,\kappa}^{-1}\left((\hat{D}\varepsilon)^{-\left(\frac{r}{q-1}+a_0\right)^{-1}}\right)$ we find that $t\geqslant\{2L,1,\tilde{T}\}$ and that the factor multiplying $q^j$ in \eqref{final lb U exponential growth} is positive, so, letting $j\to \infty$ in \eqref{final lb U exponential growth} we see that the lower bound for $U(t)$ is not finite. Hence, we proved that $U$ blows up in finite time and, as byproduct of the iteration procedure, we got the upper bound estimate for the lifespan in \eqref{upper bound for the lifespan exponential growth}. This complete the proof of Theorem \ref{Thm exponential growth}.
 
 \begin{remark}
 In the proof of \eqref{lb estimate U step j exponential growth} the assumption $r>r_{\mathrm{crit}}(b,m^2,\beta,p)$ allows to define properly $\ell_0$, since $r+\alpha_{1/2}>0$. However, the crucial point in the previous iteration argument where this assumption on the range for $r$ is used is in the representation \eqref{representation aj} for $a_j$. Indeed, the term $a_j$ allows to get a growth of exponential type in the lower bound estimates \eqref{lb estimate U step j exponential growth}. In the next subsections, we consider the limit case $r=r_{\mathrm{crit}}(b,m^2,\beta,p)$ for which the previous argument does no longer hold. A first step will be to introduce a new time-dependent functional related to $U$ and the relative iteration frame. In this new iteration frame we have to deal with just one or no exponential multiplier depending on whether we consider the case $b^2>4m^2$ or the case $b^2=4m^2$. Therefore, in Theorems \ref{Thm polynomial growth} and \ref{Thm logarithmic growth} a significant role will be played by the power $\kappa$ for the polynomial term.
 \end{remark}
 
\begin{remark} \label{Remark initial data exponential growth}

For $b>0$ and $m^2\in \left[0,\tfrac{b^2}{4}\right]$ we may weaken the sign assumptions on $u_0,u_1$ in Theorem \ref{Thm exponential growth}. In fact, it is sufficient to suppose that $\int_{\mathbb{R}^n}u_0(x)\, \mathrm{d}x, \int_{\mathbb{R}^n}u_1(x)\, \mathrm{d}x$ are nonnegative and that at least one between them is strictly positive. Indeed, under these assumptions \eqref{1st lb U} keeps to be fulfilled for $b\neq 0$ and this suffices to start the iteration argument.

On the other hand, it is interesting to consider for $b=m^2=0$ the case in which the second Cauchy data satisfy $\int_{\mathbb{R}^n}u_1(x)\, \mathrm{d}x=0$ (and, of course, $\int_{\mathbb{R}^n}u_0(x)\, \mathrm{d}x>0$). Then, Theorem \ref{Thm exponential growth} is still valid in the case $b=m^2=0$, however, the lifespan estimate in this case is the same one as in the case $b^2>4m^2$ in \eqref{upper bound for the lifespan exponential growth}. This worsening in the upper bound is caused by the fact that the lower bound for $U$ in this case is given by $U(t)\geqslant K_0 \varepsilon \, \mathrm{e}^{-\left(\frac{b}{2}-\frac{1}{2}\sqrt{b^2-4m^2}\right)t}$, i.e., without any additional linearly increasing $t$-factor on the right-hand side differently from \eqref{1st lb U}.

\end{remark}

\subsection{Case with polynomial growth: proof of Theorem \ref{Thm polynomial growth}} \label{Subsection polynomial growth}

In the present subsection, we provide the proof of Theorem \ref{Thm polynomial growth}. In this framework, the time-dependent factor $\Gamma$ in \eqref{nonlinearity} is given by \eqref{Gamma factor polynomial growth} with $\kappa>\kappa_{\mathrm{crit}}(b,m^2,\beta,p)$.

Let us multiply both sides of \eqref{Iteration frame} by $\mathrm{e}^{\alpha_1 t}$. Then, introducing the functional
\begin{align}\label{def mathcal U}
\mathcal{U}(t) \doteq \mathrm{e}^{\alpha_1 t} U(t) \qquad \mbox{for} \ t\in[0,T),
\end{align} we obtain 
\begin{align*}
\mathcal{U}(t)\geqslant C \mathrm{e}^{(\alpha_1-\alpha_2)t} \int_0^t \mathrm{e}^{(\alpha_2-\alpha_1)s}  \int_0^s \mathrm{e}^{-\alpha_1((\beta+1)p-1) \tau} \, \Gamma(\tau) (\mathcal{U}(\tau))^{(\beta+1)p} \mathrm{d}\tau \, \mathrm{d}s
\end{align*} for $t\geqslant 0$.
Differently from the previous subsection, where the role of $\alpha_{1}$ and $\alpha_2$ are interchangeable, we need to set specific values for $\alpha_1$ and $\alpha_2$. From the previous inequality it is clear that it would be beneficial to fix $\alpha_1$ in a such a way that in the $\tau$-integral the exponential factor $\mathrm{e}^{-\alpha_1((\beta+1)p-1) \tau} $ is balanced by the exponential factor $\mathrm{e}^{r_{\mathrm{crit}}(b,m^2,\beta,p)}\tau$ in $\Gamma(\tau)$. Therefore, hereafter we set
$
\alpha_1=\frac{b}{2}-\frac{1}{2}\sqrt{b^2-4m^2} \quad \mbox{and} \quad \alpha_2=\frac{b}{2}+\frac{1}{2}\sqrt{b^2-4m^2}.
$
In particular, with this choice we obtain from the previous inequality the following iteration frame for $\mathcal{U}$
\begin{align}\label{Iteration frame mathcal U}
\mathcal{U}(t)\geqslant C \mu \, \mathrm{e}^{(\alpha_1-\alpha_2)t} \int_0^t \mathrm{e}^{(\alpha_2-\alpha_1)s}  \int_0^s (1+\tau)^\kappa (\mathcal{U}(\tau))^{q} \mathrm{d}\tau \, \mathrm{d}s
\end{align} for $t\geqslant 0$, where $q=(\beta+1)p$ as in the previous subsection. Notice that the coefficient $\alpha_2-\alpha_1$ in the exponential multiplier in the $s$-integral is positive, due to our choice.

\begin{remark} Considering alternatively the functional $\widetilde{\mathcal{U}}(t)\doteq  \mathrm{e}^{\alpha_2 t}U(t)$ and switching the values of $\alpha_1$ and $\alpha_2$ with respect to the values we have just fixed, we would have found the iteration frame 
\begin{align*}
\widetilde{\mathcal{U}}(t)\geqslant C \mu  \int_0^t \mathrm{e}^{(\alpha_2-\alpha_1)s}  \int_0^s \mathrm{e}^{(\alpha_1-\alpha_2)\tau} (1+\tau)^\kappa (\widetilde{\mathcal{U}}(\tau))^{q} \mathrm{d}\tau \, \mathrm{d}s
\end{align*} for $t\geqslant 0$. Even though the structure of this iteration frame would require somehow different computations in the induction step (since the slicing procedure has to be carried out in the $\tau$-integral rather than in the $s$-integral as we will do in the next steps of the proof), the final outcome, meaning the blow-up of $u$ and the upper bound estimate for the lifespan, is exactly the same. In this sense, we can still say that the role of $\alpha_{1}$ and $\alpha_2$ are interchangeable even in this limit case for $r$.
\end{remark}

From \eqref{1st lb U} and \eqref{def mathcal U}, we get immediately the first lower bound estimates for $\mathcal{U}$, namely,
\begin{align} \label{1st lb mathcal U}
\mathcal{U}(t)\geqslant \begin{cases} K_0 \varepsilon & \mbox{if} \ b^2>4m^2, \\
K_0 \varepsilon \, (1+t) & \mbox{if} \ b^2=4m^2, 
\end{cases} 
\end{align} for $t\geqslant 0$.

From \eqref{Iteration frame mathcal U} it is clear that when $\alpha_1=\alpha_2$, that is for $b^2=4m^2$, the iteration procedure which we use to establish the sequence of lower bound estimates for $\mathcal{U}$ is quite different. Indeed, depending on whether or not an unbounded exponential multiplier is present in the $s$-integral we might need to apply the slicing procedure or not. Hence, we will consider separately the cases $b^2>4m^2$ and $b^2=4m^2$.

\subsubsection{Case with polynomial growth: sub-case with dominant damping} 

In this case $\alpha_1\neq \alpha_2$ so that $\alpha_2-\alpha_1=\sqrt{b^2-4m^2}>0$. Since in the iteration frame for $\mathcal{U}$ given by \eqref{Iteration frame mathcal U} we have the exponential multiplier $\mathrm{e}^{(\alpha_2-\alpha_1)s}$, we have to modify the choice of the parameters $\{L_j\}_{j\in\mathbb{N}}$ characterizing the slicing procedure with respect to Subsection \ref{Subsection exponential growth}. Formally, $L_j$ is defined as in \eqref{def Lj exponential growth}, however, the coefficients $\{\ell_k\}_{k\in\mathbb{N}}$ are given in this case by
\begin{align*}
\ell_0 & \doteq (\alpha_2-\alpha_1)^{-1}, \\
\ell_k & \doteq 1+ q^{-k} \qquad \qquad \mbox{for any} \ k\geqslant 1 .
\end{align*} Since $\ell_k>1$ for any $k\in\mathbb{N}\setminus \{0\}$, also in this case we have $L_j \uparrow$. Moreover, we keep using the notation $L\doteq \lim_{j\to\infty} L_j$ (the convergence of this infinite product can be proved exactly as in the previous case).

Let us prove now the following sequence of lower bound estimates for $\mathcal{U}$:
\begin{align}\label{lower bound estimates mathcal U j}
\mathcal{U}(t)\geqslant C_j (t-L_j)^{b_j} (1+t)^{-\beta_j} \qquad \mbox{for} \ t\geqslant L_j \ \mbox{and for any} \ j\in\mathbb{N},
\end{align} where $\{C_j\}_{j\in\mathbb{N}}$, $\{b_j\}_{j\in\mathbb{N}}$, $\{\beta_j\}_{j\in\mathbb{N}}$ are sequences of nonnegative real numbers to be determined during the inductive argument. From \eqref{1st lb mathcal U}, we have that \eqref{lower bound estimates mathcal U j} for $j=0$, provided that $C_0\doteq K_0\varepsilon$, $b_0\doteq 0$ and $\beta_0\doteq 0$.

Let us prove the induction step. Plugging \eqref{lower bound estimates mathcal U j} in \eqref{Iteration frame mathcal U} we have
\begin{align*}
\mathcal{U}(t) & \geqslant \mu  C \mathrm{e}^{(\alpha_1-\alpha_2)t} \int_{L_j}^t \mathrm{e}^{(\alpha_2-\alpha_1)s}  \int_{L_j}^s (1+\tau)^\kappa (\mathcal{U}(\tau))^{q} \mathrm{d}\tau \, \mathrm{d}s \\
& \geqslant \mu  C C_j^q (1+t)^{-(\kappa_-+q\beta_j)} \,\mathrm{e}^{(\alpha_1-\alpha_2)t} \int_{L_j}^t \mathrm{e}^{(\alpha_2-\alpha_1)s}  \int_{L_j}^s  (\tau-L_j)^{\kappa_++b_j q} \mathrm{d}\tau \, \mathrm{d}s \\
& = \mu  C C_j^q (1+\kappa_++q b_j)^{-1}(1+t)^{-(\kappa_-+q\beta_j)} \,\mathrm{e}^{(\alpha_1-\alpha_2)t} \int_{L_j}^t \mathrm{e}^{(\alpha_2-\alpha_1)s}   (s-L_j)^{1+\kappa_++qb_j} \, \mathrm{d}s
\end{align*} for $t\geqslant L_j$, where $\kappa_+,\kappa_-$ denote the positive and the negative part of $\kappa$, respectively. For $t\geqslant L_{j+1}$, it is possible to shrink the domain of integration to $[t/\ell_{j+1},t]$ in the last integral, obtaining
\begin{align*}
\mathcal{U}(t) & \geqslant \mu  C C_j^q (1+\kappa_++q b_j )^{-1}(1+t)^{-(\kappa_-+q\beta_j)} \,\mathrm{e}^{(\alpha_1-\alpha_2)t} \int_{\tfrac{L_j t}{L_{j+1}}}^t \mathrm{e}^{(\alpha_2-\alpha_1)s}   (s-L_j)^{1+\kappa_++qb_j} \, \mathrm{d}s \\
& \geqslant \frac{\mu  C C_j^q}{(1+\kappa_++q b_j ) \, \ell_{j+1}^{1+\kappa_++qb_j}}  (t-L_j)^{1+\kappa_++qb_j } (1+t)^{-(\kappa_-+q\beta_j)} \,\mathrm{e}^{(\alpha_1-\alpha_2)t} \int_{\tfrac{L_j t}{L_{j+1}}}^t \mathrm{e}^{(\alpha_2-\alpha_1)s}  \, \mathrm{d}s \\
& \geqslant \frac{\mu  C (\alpha_2-\alpha_1)^{-1} C_j^q}{(1+\kappa_++q b_j ) \, \ell_{j+1}^{1+\kappa_++qb_j}}  (t-L_j)^{1+\kappa_++qb_j } (1+t)^{-(\kappa_-+q\beta_j)} \left(1- \mathrm{e}^{-(\alpha_2-\alpha_1)\big(1-\frac{1}{\ell_{j+1}}\big)t} \right).
\end{align*} Using the estimate 
\begin{align}
1- \mathrm{e}^{-(\alpha_2-\alpha_1)\big(1-\frac{1}{\ell_{j+1}}\big)t} & \geqslant 1- \mathrm{e}^{-(\alpha_2-\alpha_1)(\ell_{j+1}-1)L_{j}} 
 \geqslant 1- \mathrm{e}^{-(\alpha_2-\alpha_1)(\ell_{j+1}-1)L_{0}} = 1- \mathrm{e}^{-(\ell_{j+1}-1)}\notag \\
 &\geqslant (\ell_{j+1}-1)\left(1-\tfrac{1}{2}(\ell_{j+1}-1) \right) \geqslant (q-\tfrac 12)q^{-2(j+1)} \label{lb exponential factor}
\end{align}
for $t\geqslant L_{j+1}$, we find
\begin{align*}
\mathcal{U}(t) &  \geqslant \frac{\mu  C (q-\tfrac 12) C_j^q}{(\alpha_2-\alpha_1)(1+\kappa_++q b_j ) \, \ell_{j+1}^{1+\kappa_++qb_j} q^{2(j+1)}}  (t-L_j)^{1+\kappa_++qb_j } (1+t)^{-(\kappa_-+q\beta_j)},
\end{align*} which is exactly \eqref{lower bound estimates mathcal U j} for $j+1$, provided that
\begin{align}
C_{j+1} & \doteq \frac{\mu  C (q-\tfrac 12)}{(\alpha_2-\alpha_1)} (1+\kappa_++q b_j )^{-1} \, \ell_{j+1}^{-(1+\kappa_++qb_j)} q^{-2(j+1)}C_j^q  \label{def Cj polynomial growth},\\
b_{j+1} & \doteq 1+\kappa_++qb_j , \quad \beta_{j+1}\doteq \kappa_-+q\beta_j. \label{def bj beta j polynomial growth}
\end{align} By employing recursively \eqref{def bj beta j polynomial growth} among two consecutive terms from the sequences $\{b_j\}_{j\in\mathbb{N}}$ and $\{\beta_j\}_{j\in\mathbb{N}}$, we get
\begin{align}
b_j & = \frac{1+\kappa_+}{q-1}q^j - \frac{1+\kappa_+}{q-1}, \label{representation bj polynomial growth} \\
\beta_j & = \frac{\kappa_-}{q-1} q^j - \frac{\kappa_-}{q-1}, \label{representation beta j polynomial growth}
\end{align} where we used $b_0=\beta_0=0$.
Thanks to \eqref{def bj beta j polynomial growth} and \eqref{representation bj polynomial growth}, we have
\begin{align}\label{upper bound b j+1 poly growth}
1+\kappa_++q b_j = b_{j+1}\leqslant \frac{1+\kappa_+}{q-1} \, q^{j+1} .
\end{align} Moreover,
\begin{align*}
\lim_{j\to\infty} \ell_{j+1}^{b_{j+1}}= \lim_{j\to \infty} \exp(b_{j+1}\ln \ell_{j+1}) = \lim_{j\to \infty} \exp\left(\tfrac{1+\kappa_+}{q-1}q^{j+1}\ln\left(1+q^{-(j+1)}\right)\right)= \exp\left(\tfrac{1+\kappa_+}{q-1}\right)
\end{align*} implies the existence of a constant $M_2=M_2(\beta,p,\kappa)$ such that $ \ell_{j+1}^{b_{j+1}}\leqslant M_2$ for any $j\in\mathbb{N}$. Combining this last uniform upper bound with \eqref{def Cj polynomial growth} and \eqref{upper bound b j+1 poly growth}, we obtain
\begin{align*}
C_{j+1} & \doteq \underbrace{\frac{\mu  C (q-\tfrac 12)(q-1)}{(\alpha_2-\alpha_1)(1+\kappa_+) M_2}}_{B\doteq }  q^{-3(j+1)}C_j^q 
\end{align*} for any $j\in\mathbb{N}$.
Applying the logarithmic function to both sides of the  inequality  $C_j\geqslant B q^{-3j} C_{j-1}^q$ and, then, using iteratively the resulting inequality, we find 
\begin{align*}
\ln C_j &  \geqslant q\ln C_{j-1}-3j \ln q+\ln b \geqslant q^2\ln C_{j-1}-3(j+(j-1)q) \ln q+ (1+q)\ln B \\
 &  \geqslant \ldots \geqslant q^j \ln C_0-3 \left(\,\sum_{k=0}^{j-1}(j-k)q^k\right)\ln q+\left(\,\sum_{k=0}^{j-1}q^k\right)\ln B \\
  & \geqslant q^j\left( \ln C_0-\frac{3q \ln q}{(q-1)^2}+\frac{\ln B}{q-1}\right)+\frac{3q \ln q}{(q-1)^2}
+\frac{3\ln q}{q-1} j-\frac{\ln B}{q-1},
\end{align*} where in the last step we used \eqref{summation identity}.

 Let $j_1=j_1(n,c,H,b,m^2,\mu,\kappa,\beta,p,R)\in\mathbb{N}$ be the smallest integer such that $j_1\geqslant \frac{\ln B}{3\ln q}-\frac{q}{q-1}$. Hence, for any $j\geqslant j_1$ it holds
\begin{align}\label{exponential lb for Cj poly growth}
\ln C_j &  \geqslant q^j\left( \ln (K_0 \varepsilon)-\frac{3q \ln q}{(q-1)^2}+\frac{\ln B}{q-1}\right)= q^j \ln (\tilde{B}\varepsilon),
\end{align} where $\tilde{B}\doteq K_0 q^{-3q/(q-1)^2}B^{1/(q-1)}$.
Since $L_j\uparrow L$ as $j\to \infty$, in particular \eqref{lower bound estimates mathcal U j} is true for $t\geqslant L$ and any $j\in\mathbb{N}$. Combining \eqref{lower bound estimates mathcal U j}, \eqref{representation bj polynomial growth}, \eqref{representation beta j polynomial growth} and \eqref{exponential lb for Cj poly growth}, for $t\geqslant L$ and $j\geqslant j_1$ we arrive at
\begin{align*}
\mathcal{U}(t) & \geqslant \exp\left(q^j \left(\ln (\tilde{B}\varepsilon)+\tfrac{1+\kappa_+}{q-1}\ln (t-L)-\tfrac{\kappa_-}{q-1}\ln (1+t) \right)\right) (t-L)^{-\frac{1+\kappa_+}{q-1}} (1+t)^{-\frac{\kappa_-}{q-1}}. 
\end{align*} Then, for $t\geqslant \max\{2L,1\}$ and $j\geqslant j_1$, by using $\kappa=\kappa_+-\kappa_-$, from the previous estimate we have
\begin{align} \label{final lb U poly growth}
\mathcal{U}(t) & \geqslant \exp\left(q^j \ln \left(\hat{B}\varepsilon t^{\frac{1+\kappa}{q-1}}\right)\right)  (t-L)^{-\frac{1+\kappa_+}{q-1}} (1+t)^{-\frac{\kappa_-}{q-1}},
\end{align} where $\hat{B}\doteq 2^{-(1+\kappa_++\kappa_-)/(q-1) }\tilde{B}$.

 The logarithmic factor multiplying $q^j$ in \eqref{final lb U poly growth} is strictly positive if and only if $ t>(\hat{B}\varepsilon)^{-\frac{q-1}{1+\kappa}}$. We set $\varepsilon_0=\varepsilon_0(n,c,H,b,m^2,\beta,p,\mu,\kappa,u_0,u_1,R)~>~0$ sufficiently small so that $$(\hat{B}\varepsilon_0)^{-\frac{q-1}{1+\kappa}} \geqslant \max\{2L,1\}.$$ Then, for any $\varepsilon\in(0,\varepsilon_0]$ and any $t>(\hat{B}\varepsilon)^{-\frac{q-1}{1+\kappa}}$ we obtain that $t\geqslant\{2L,1\}$ and that the factor multiplying $q^j$ in \eqref{final lb U poly growth} is positive, thus, taking the limit as $j\to \infty$ in \eqref{final lb U poly growth} we have that the lower bound for $\mathcal{U}(t)$ is not finite. Hence, we proved that $\mathcal{U}$ blows up in finite time as well and, as byproduct of the former iteration procedure, the upper bound estimate for the lifespan in \eqref{upper bound for the lifespan polynomial growth} has been proven when $b^2>4m^2$.

\subsubsection{Case with polynomial growth: sub-case with balanced damping and mass} 

When $b^2=4m^2$ in \eqref{Iteration frame mathcal U} the exponential terms disappear.

Indeed, in this special case the iteration frame for $\mathcal{U}$ is given by  
\begin{align} \label{Iteration frame mathcal U balanced case}
\mathcal{U}(t)\geqslant C \mu  \int_0^t   \int_0^s (1+\tau)^\kappa (\mathcal{U}(\tau))^{q} \mathrm{d}\tau \, \mathrm{d}s
\end{align} for $t\geqslant 0$. Combining \eqref{Iteration frame mathcal U balanced case} with the lower bound estimate for $\mathcal{U}$ in \eqref{1st lb mathcal U}, from \cite[Lemma 1]{LaiSchTak20}  we see that $\mathcal{U}$ blows up in finite time provided that $\gamma\doteq 2[(q-1)+\kappa+2]>0$, that is, for $\kappa >-1-q =\kappa_{\mathrm{crit}}(b,m^2,\beta,p)$. Moreover, the upper bound estimate $T(\varepsilon)\lesssim \widetilde{T}(\varepsilon)$ holds for the lifespan, where $\widetilde{T}$ is defined through the relation $\varepsilon \widetilde{T}^{\frac{\gamma}{2(q-1)}}=1$. Consequently, $T(\varepsilon)\lesssim \varepsilon^{-(\frac{\kappa+2}{q-1}+1)^{-1}}$, which is exactly the upper bound estimate in \eqref{upper bound for the lifespan polynomial growth} for the case $b^2=4m^2$.

\begin{remark} \label{Remark initial data polynomial growth}
We  emphasize that for $b>0$ and $m^2\in \left[0,\tfrac{b^2}{4}\right]$ in order to prove Theorem \ref{Thm polynomial growth}, concerning the sign assumptions for the Cauchy data it is sufficient to require that $\int_{\mathbb{R}^n}u_0(x)\, \mathrm{d}x, \int_{\mathbb{R}^n}u_1(x)\, \mathrm{d}x$ are nonnegative and that at least one between them is strictly positive (analogously to what we pointed out in Remark \ref{Remark initial data exponential growth} for Theorem \ref{Thm exponential growth}). Indeed, under these assumptions \eqref{1st lb mathcal U} still holds true if $b>0$.

On the other hand, if we have $\int_{\mathbb{R}^n}u_1(x)\, \mathrm{d}x=0$ (and, of course, $\int_{\mathbb{R}^n}u_0(x)\, \mathrm{d}x>0$), Theorem \ref{Thm polynomial growth} still holds in the case $b=m^2=0$, however, the range for $\kappa$  becomes $\kappa>-2$ and the lifespan estimate is in this case changes to $$T(\varepsilon)\lesssim \varepsilon^{-\frac{q-1}{\kappa+2}},$$ due to the fact that the lower bound for $\mathcal{U}$ in this case is given by $\mathcal{U}(t)\geqslant K_0 \varepsilon$ (i.e., without any linear increasing $t$-factor on the right-hand side).
\end{remark}

\subsection{Case with logarithmic growth: proof of Theorem \ref{Thm logarithmic growth}} \label{Subsection log growth}

In this section, we prove Theorem \ref{Thm logarithmic growth} by showing that the functional $\mathcal{U}$ introduced in Subsection \ref{Subsection polynomial growth} blows up even when the $\Gamma$ factor in \eqref{nonlinearity} is given by \eqref{Gamma factor logarithmic growth} with threshold values both for the exponential factor and for the polynomial factor. The iteration frame is the one given in \eqref{Iteration frame mathcal U}. Nonetheless, we will employ it for deriving different kinds of lower bound estimates for $\mathcal{U}$, depending on whether we work with $b^2>4m^2$ or with $b^2=4m^2$. We emphasize that in this final case we still need to apply a slicing procedure to handle logarithmic terms in the $s$-integral. In the dominant damping case (i.e., for $b^2>4m^2$) the slicing procedure will enable us to control both the logarithmic factors and the exponential multiplier, while in the balanced case $b^2=4m^2$ the exponential term disappears, so the slicing procedure will serve to deal with the logarithmic terms only.

In the first case we work with the same sequence $\{L_j\}_{j\in\mathbb{N}}$ as the one defined in Subsection \ref{Subsection polynomial growth}. Whilst in the case $b^2=4m^2$ we consider a simpler sequence $\{L_j\}_{j\in\mathbb{N}}$ which is analogous to the one introduced for the first time in \cite{AKT00}.

\subsubsection{Case with logarithmic growth: sub-case with dominant damping}

When the damping term is dominant ($b^2>4m^2$) and for $\kappa=\kappa_{\mathrm{crit}}(b,m^2,\beta,p)=-1$, the iteration frame in \eqref{Iteration frame mathcal U} can be rewritten as follows:
\begin{align}\label{Iteration frame mathcal U log growth dominant damp}
\mathcal{U}(t)\geqslant \mu C  \,\mathrm{e}^{(\alpha_1-\alpha_2)t} \int_0^t \mathrm{e}^{(\alpha_2-\alpha_1)s}  \int_0^s (1+\tau)^{-1} (\mathcal{U}(\tau))^{q} \mathrm{d}\tau \, \mathrm{d}s \qquad \mbox{for} \ t\geqslant 0. 
\end{align}
The next step is to show the following sequence of lower bound estimates for $\mathcal{U}$
\begin{align} \label{lb estimates mathcal U j log growth dominant damp}
\mathcal{U}(t)\geqslant C_j \left(\ln\left(\frac{t}{L_{j}}\right)\right)^{d_j} \qquad \mbox{for} \  t\geqslant L_j \ \mbox{and any} \ j\in\mathbb{N}, 
\end{align} where $\{C_j\}_{j\in\mathbb{N}}$, $\{d_j\}_{j\in\mathbb{N}}$ are sequences of nonnegative real numbers to be determined iteratively, and $\{L_j\}_{j\in\mathbb{N}}$ is defined as in Subsection \ref{Subsection polynomial growth} in the case $b^2>4m^2$.

Clearly, from \eqref{1st lb mathcal U} we get \eqref{lb estimates mathcal U j log growth dominant damp} with $C_0\doteq K_0 \varepsilon$ and $d_0\doteq 0$. Let us proceed now with the induction step. Plugging \eqref{lb estimates mathcal U j log growth dominant damp} in \eqref{Iteration frame mathcal U log growth dominant damp}, for $t\geqslant L_j$ we obtain
\begin{align*}
\mathcal{U}(t) & \geqslant  \mu  C \, \mathrm{e}^{(\alpha_1-\alpha_2)t} \int_{L_j}^t \mathrm{e}^{(\alpha_2-\alpha_1)s}  \int_{L_j}^s (1+\tau)^{-1} (\mathcal{U}(\tau))^{q} \mathrm{d}\tau \, \mathrm{d}s \\
 & \geqslant  \mu  C C_j^q \left(1+\frac{1}{L_j}\right)^{-1}\mathrm{e}^{(\alpha_1-\alpha_2)t} \int_{L_j}^t \mathrm{e}^{(\alpha_2-\alpha_1)s}  \int_{L_j}^s \tau^{-1}\left(\ln\left(\frac{\tau}{L_{j}}\right)\right)^{q d_j} \mathrm{d}\tau \, \mathrm{d}s \\
 & \geqslant  \mu  C C_j^q \left(1+\frac{1}{\ell_0}\right)^{-1} (1+qd_j)^{-1}\mathrm{e}^{(\alpha_1-\alpha_2)t} \int_{L_j}^t \mathrm{e}^{(\alpha_2-\alpha_1)s}  \left(\ln\left(\frac{s}{L_{j}}\right)\right)^{1+q d_j} \, \mathrm{d}s,
\end{align*} where in the second inequality we used $1+\tau\leqslant (1+L_j^{-1})\tau$ for $\tau\geqslant L_j$. For $t\geqslant L_{j+1}$, by cutting away a slice from the domain of integration, we have
\begin{align*}
\mathcal{U}(t) & \geqslant  \mu  C C_j^q \left(1+\frac{1}{\ell_0}\right)^{-1} (1+qd_j)^{-1}\mathrm{e}^{(\alpha_1-\alpha_2)t} \int_{\tfrac{L_j t}{L_{j+1}}}^t \mathrm{e}^{(\alpha_2-\alpha_1)s}  \left(\ln\left(\frac{s}{L_{j}}\right)\right)^{1+q d_j} \, \mathrm{d}s \\
& \geqslant  \mu  C C_j^q \left(1+\frac{1}{\ell_0}\right)^{-1} (1+qd_j)^{-1} \left(\ln\left(\frac{t}{L_{j+1}}\right)\right)^{1+q d_j} \mathrm{e}^{(\alpha_1-\alpha_2)t} \int_{\tfrac{L_j t}{L_{j+1}}}^t \mathrm{e}^{(\alpha_2-\alpha_1)s}  \, \mathrm{d}s \\
& = \mu  C (\alpha_2-\alpha_1)^{-2}\left(\ell_0+1\right)^{-1} C_j^q  (1+qd_j)^{-1} \left(\ln\left(\frac{t}{L_{j+1}}\right)\right)^{1+q d_j}\! \left(1-\mathrm{e}^{-(\alpha_2-\alpha_1)\big(1-\frac{1}{\ell_{j+1}}\big)t} \right) \\
& \geqslant \mu  C \left(q-\tfrac 12\right)(\alpha_2-\alpha_1)^{-2}\left(\ell_0+1\right)^{-1} C_j^q  (1+qd_j)^{-1} q^{-2(j+1)}\left(\ln\left(\frac{t}{L_{j+1}}\right)\right)^{1+q d_j} ,
\end{align*} where in the last inequality we used \eqref{lb exponential factor}. The previous chain of inequalities provides exactly \eqref{lb estimates mathcal U j log growth dominant damp} for $j+1$ by setting
\begin{align}
C_{j+1} & \doteq \frac{\mu  C \left(q-\tfrac 12\right) C_j^q }{(\alpha_2-\alpha_1)^{2}\left(\ell_0+1\right)(1+qd_j) \, q^{2(j+1)}},  \label{def C j log growth dominant damp} \\
d_{j+1} &\doteq 1+q d_j. \label{def dj log growth dominant damp}
\end{align} Analogously to what we have done in the previous subsections, we derive first an explicit representation for $d_{j}$ and then we determine a suitable lower bound for $C_j$ when $j$ is large enough. By using recursively \eqref{def dj log growth dominant damp} and $d_0=0$, we get
\begin{align}\label{representation dj dominant damp}
d_j =  1+q d_{j-1} = \sum_{k=0}^{j-1} q^k + d_0 q^j = \frac{q^j-1}{q-1}.
\end{align} Therefore, \eqref{def dj log growth dominant damp} and \eqref{representation dj dominant damp} imply that $(1+qd_j)^{-1}\geqslant (q-1) q^{-(j+1)}$. Consequently,
\begin{align*}
C_{j+1} & \doteq \underbrace{\frac{\mu  C \left(q-\tfrac 12\right) (q-1) }{(\alpha_2-\alpha_1)^{2}\left(\ell_0+1\right)}}_{E\doteq } q^{-3(j+1)} C_j^q. 
\end{align*} Similarly as we did in Subsection \ref{Subsection polynomial growth}, from the inequality $C_j\geqslant E q^{-3j} C_{j-1}$ we get
\begin{align*}
\ln C_j &  \geqslant q^j\left( \ln C_0-\frac{3q \ln q}{(q-1)^2}+\frac{\ln E}{q-1}\right)+\frac{3q \ln q}{(q-1)^2}
+\frac{3\ln q}{q-1} j-\frac{\ln E}{q-1}.
\end{align*}
Setting $j_2=j_2(n,c,H,b,m^2,\mu,\beta,p,R)\in\mathbb{N}$ to be the smallest integer such that $j_2\geqslant \frac{\ln E}{3\ln q}-\frac{q}{q-1}$, for any $j\geqslant j_2$ we get
\begin{align}\label{exponential lb for Cj log growth}
\ln C_j &  \geqslant q^j\left( \ln (K_0 \varepsilon)-\frac{3q \ln q}{(q-1)^2}+\frac{\ln E}{q-1}\right)= q^j \ln (\tilde{E}\varepsilon),
\end{align} where $\tilde{E}\doteq K_0 q^{-3q/(q-1)^2}E^{1/(q-1)}$.
Combining \eqref{lb estimates mathcal U j log growth dominant damp}, \eqref{representation dj dominant damp} and \eqref{exponential lb for Cj log growth}, for $t\geqslant L =\lim_{j\to \infty}L_{j}$ and $j\geqslant j_2$ we find
\begin{align}
\mathcal{U}(t) & \geqslant \exp\left( q^j \ln (\tilde{E}\varepsilon)\right) \left(\ln\left(\frac{t}{L}\right)\!\right)^{ \frac{q^j-1}{q-1}} \notag \\ 
&= \exp\left( q^j \left(\ln (\tilde{E}\varepsilon)+\frac{1}{q-1}\ln\left(\ln\left(\frac{t}{L}\right)\!\right)\!\right)\right) \left(\ln\left(\frac{t}{L}\right)\!\right)^{ -\frac{1}{q-1}} \notag \\
&= \exp\left( q^j \ln \left(\tilde{E}\varepsilon \left(\ln\left(\frac{t}{L}\right)\!\right)^{\frac{1}{q-1}}\right)\right)\left(\ln\left(\frac{t}{L}\right)\!\right)^{ -\frac{1}{q-1}}. \label{final lb mathcal U log growth}
\end{align} 
We remark that the logarithmic factor multiplying $q^j$ in \eqref{final lb mathcal U log growth} is strictly positive if and only if $t>L\exp(\tilde{E}\varepsilon)^{-(q-1)}$.

Now we fix $\varepsilon_0=\varepsilon_0(n,c,H,b,m^2,\beta,p,\mu,u_0,u_1,R)$ such that $\exp(\tilde{E}\varepsilon_0)^{-(q-1)}>1$. Then, for any $\varepsilon\in(0,\varepsilon_0]$ and any $t >L\exp(\tilde{E}\varepsilon)^{-(q-1)}$ the right-hand side of \eqref{final lb mathcal U log growth} diverges as $j\to \infty$, so $\mathcal{U}(t)$ cannot be finite. Summarizing,  we proved the blow-up of $\mathcal{U}$ in finite time and the upper bound estimate in \eqref{upper bound for the lifespan logarithmic growth} when $b^2>4m^2$.

\subsubsection{Case with logarithmic growth: sub-case with balanced damping and mass} 

If $b^2=4m^2$ and $\kappa=\kappa_{\mathrm{crit}}(b,m^2,\beta,p)=-1- q$, we may rewrite the iteration frame  \eqref{Iteration frame mathcal U} as follows:
\begin{align}\label{Iteration frame mathcal U log growth balanced}
\mathcal{U}(t)\geqslant \mu C   \int_0^t  \int_0^s (1+\tau)^{-1-q} (\mathcal{U}(\tau))^{q} \mathrm{d}\tau \, \mathrm{d}s \qquad \mbox{for} \ t\geqslant 0. 
\end{align}
The next step is to show the following sequence of lower bound estimates for $\mathcal{U}$
\begin{align} \label{lb estimates mathcal U j log growth balanced}
\mathcal{U}(t)\geqslant C_j t \left(\ln\left(\frac{t}{L_{j}}\right)\!\right)^{d_j} \qquad \mbox{for} \  t\geqslant L_j \ \mbox{and any} \ j\in\mathbb{N}, 
\end{align} where $\{C_j\}_{j\in\mathbb{N}}$, $\{d_j\}_{j\in\mathbb{N}}$ are suitable sequences of nonnegative real numbers, and 
\begin{align}\label{Lj slicing procedure log growth balanced}
 L_j \doteq 
 2-2^{-j} \quad \mbox{for} \  j\in\mathbb{N}.
 \end{align}

\begin{remark} The choice of the sequence $\{L_j\}_{j\in\mathbb{N}}$ in \eqref{Lj slicing procedure log growth balanced} is done in order to handle the logarithmic factors in the $s$-integral. Notice that in this case  no exponential multiplier appears, so the construction of the parameters characterizing the slicing procedure is simpler than in the previous proofs and it is inspired by the one from \cite[Section 6]{AKT00}.
\end{remark}

We begin by observing that \eqref{lb estimates mathcal U j log growth balanced} for $j=0$ follows from \eqref{1st lb mathcal U} with $C_0\doteq K_0 \varepsilon$ and $d_0\doteq 0$. Let us proceed with the inductive step. If we plug \eqref{lb estimates mathcal U j log growth balanced} in \eqref{Iteration frame mathcal U log growth balanced}, for $t\geqslant L_j$ it results
\begin{align*}
\mathcal{U}(t) & \geqslant \mu C \left(1+\frac{1}{L_j}\right)^{-1-q}  \int_{L_j}^t  \int_{L_j}^s \tau^{-1-q} (\mathcal{U}(\tau))^{q} \mathrm{d}\tau \, \mathrm{d}s \\
& \geqslant \mu C C_j^q \left(1+\frac{1}{L_j}\right)^{-1-q}  \int_{L_j}^t  \int_{L_j}^s \tau^{-1} \left(\ln\left(\frac{\tau}{L_{j}}\right)\!\right)^{q d_j}  \mathrm{d}\tau \, \mathrm{d}s \\
& \geqslant \mu C C_j^q \left(1+\frac{1}{L_j}\right)^{-1-q} (1+q d_j)^{-1} \int_{L_j}^t  \left(\ln\left(\frac{s}{L_{j}}\right)\!\right)^{1+q d_j}   \mathrm{d}s.
\end{align*} Hence, using the slicing procedure, for $t\geqslant L_{j+1}$ we have
\begin{align*}
\mathcal{U}(t) & \geqslant \mu C C_j^q \left(1+\frac{1}{L_j}\right)^{-1-q} (1+q d_j)^{-1} \int_{\tfrac{L_j t}{L_{j+1}}}^t  \left(\ln\left(\frac{s}{L_{j}}\right)\!\right)^{1+q d_j}   \mathrm{d}s \\
& = \mu C C_j^q \left(1+\frac{1}{L_j}\right)^{-1-q} (1+q d_j)^{-1} \left(1-\frac{L_j}{L_{j+1}}\right) t \left(\ln\left(\frac{t}{L_{j+1}}\right)\!\right)^{1+q d_j}, 
\end{align*} which is \eqref{lb estimates mathcal U j log growth balanced} for $j+1$ provided that we set
\begin{align}
C_{j+1} & \doteq \mu C (1+q d_j)^{-1} \left(1+\frac{1}{L_j}\right)^{-1-q}  \left(1-\frac{L_j}{L_{j+1}}\right) C_j^q ,  \label{def C j log growth balanced} \\
d_{j+1} &\doteq 1+q d_j. \label{def dj log growth balanced}
\end{align} 
In a complete analogous way as in the previous case, we derive the representation \eqref{representation dj dominant damp} for $d_j$ even in this case. Hence, using 
\begin{align*}
\left(1+\frac{1}{L_j}\right)^{-1-q} \geqslant 2^{-1-q} , \quad 1-\frac{L_j}{L_{j+1}}\geqslant 2^{-(j+2)}
\end{align*} for any $j\in\mathbb{N}$ and \eqref{representation dj dominant damp}, we arrive at
\begin{align*}
C_{j+1}\geqslant \underbrace{2^{-2-q} \mu C (q-1)}_{F\doteq } (2q)^{-(j+1) } C_j^q.
\end{align*} Repeating similar computations as in the previous proofs, from the inequality $C_j\geqslant F (2q)^{-j } C_{j-1}^q$, we can derive the following lower bound for $\ln C_j$
\begin{align}\label{exponential lb for Cj log growth balanced}
\ln C_j &  \geqslant q^j\left( \ln (K_0 \varepsilon)-\frac{q \ln (2q)}{(q-1)^2}+\frac{\ln F}{q-1}\right)= q^j \ln (\tilde{F}\varepsilon),
\end{align} for any $j\geqslant j_3$, where $j_3=j_3(n,c,H,b,m^2,\mu,\beta,p,R)\in\mathbb{N}$ is the smallest integer such that $j_3\geqslant \frac{\ln F}{\ln (2q)}-\frac{q}{q-1}$ and $\tilde{F}\doteq K_0 (2q)^{-q/(q-1)^2}F^{1/(q-1)}$.

Combining \eqref{lb estimates mathcal U j log growth balanced} and \eqref{exponential lb for Cj log growth balanced}, for $t\geqslant 2 =\lim_{j\to \infty}L_{j}$ and $j\geqslant j_3$ we have
\begin{align}
\mathcal{U}(t) & \geqslant \exp\left( q^j \ln (\tilde{F}\varepsilon)\right) t \left(\ln\left(\frac{t}{2}\right)\!\right)^{ \frac{q^j-1}{q-1}} \notag \\ 
&= \exp\left( q^j \ln \left(\tilde{F}\varepsilon \left(\ln\left(\frac{t}{2}\right)\!\right)^{\frac{1}{q-1}}\right)\!\right) t \left(\ln\left(\frac{t}{2}\right)\!\right)^{ -\frac{1}{q-1}}. \label{final lb mathcal U log growth balanced}
\end{align} 
In \eqref{final lb mathcal U log growth balanced} the logarithmic factor multiplying $q^j$  is strictly positive if and only if $t>2\exp(\tilde{F}\varepsilon)^{-(q-1)}$, therefore, we fix $\varepsilon_0=\varepsilon_0(n,c,H,b,m^2,\beta,p,\mu,u_0,u_1,R)$ such that $\exp(\tilde{F}\varepsilon_0)^{-(q-1)}>1$. Thus, for any $\varepsilon\in(0,\varepsilon_0]$ and any $t >2\exp(\tilde{F}\varepsilon)^{-(q-1)}$ the right-hand side of \eqref{final lb mathcal U log growth balanced} diverges as $j\to \infty$ and, in particular, $\mathcal{U}(t)$ cannot be finite. In conclusion, we showed the blow-up of $\mathcal{U}$ in finite time and the upper bound estimate in \eqref{upper bound for the lifespan logarithmic growth} for $b^2=4m^2$.

\begin{remark} Analogously to what we pointed out in Remarks \ref{Remark initial data exponential growth} and \ref{Remark initial data polynomial growth} for Theorems \ref{Thm exponential growth} and \ref{Thm polynomial growth}, respectively, it is possible to weaken the sign assumptions on the Cauchy data in the statement of Theorem \ref{Thm logarithmic growth}. More precisely, for $b>0$ and $m^2\in \left[0,\frac{b^2}{4}\right]$ assuming that that $\int_{\mathbb{R}^n}u_0(x)\, \mathrm{d}x, \int_{\mathbb{R}^n}u_1(x)\, \mathrm{d}x$ are nonnegative and that at least one between them is strictly positive, then, the blow-up result from Theorem \ref{Thm logarithmic growth} is still valid. 

On the other hand, in the limit case $b=m^2=0$ if $ \int_{\mathbb{R}^n}u_1(x)\, \mathrm{d}x>0$, then, it holds the same blow-up result as in Theorem \ref{Thm logarithmic growth}; while if  $\int_{\mathbb{R}^n}u_1(x)\, \mathrm{d}x=0$ and $ \int_{\mathbb{R}^n}u_0(x)\, \mathrm{d}x>0$ then the result has to be modified  accordingly to Remark \ref{Remark initial data polynomial growth} with the threshold value $\kappa=-1-q$ replaced by $\kappa=-2$. 
\end{remark}

\section{Models in anti-de Sitter spacetime}

\subsection{Derivation of the iteration frame}\label{Section iteration frame anti-dS}

The derivation of the iteration frame for \eqref{anti deSitter semi} can be done in a complete analogous way as we did for \eqref{deSitter semi} in Subsection \ref{Section iteration frame}.

If $v$ is a local in time solution to \eqref{anti deSitter semi}, denoting 
\begin{align*}
V(t)\doteq \int_{\mathbb{R}^n} v(t,x) \, \mathrm{d}x \qquad \mbox{for} \ t\in(0,T),
\end{align*} then, the iteration frame involves this functional $V$.

Fixed $t\in (0,T)$ we consider a cutoff function $\varphi\in\mathcal{C}([0,T)\times\mathbb{R}^n)$ that localizes the support of $v$ on the strip $[0,t]\times \mathbb{R}^n$, that is, $\varphi=1$ on $\{(s,x)\in [0,t]\times \mathbb{R}^n:|x|\leqslant R+cH^{-1}(\mathrm{e}^{sH}-1)\}$. Consequently, employing this $\varphi$ in \eqref{def weak sol int rel anti dS} and differentiating with respect to $t$ the resulting relation, we have
\begin{align*}
V''(t)+bV'(t)+m^2V(t)= \Gamma(t)\left( \int_{\mathbb{R}^n}  |v(t,x)|^p \mathrm{d}x\right)^{\beta+1}.
\end{align*} which is formally identical to \eqref{ODE for F}. By using the same factorization of the operator $\frac{\mathrm{d}^2}{\mathrm{d}t^2}+b\frac{\mathrm{d}}{\mathrm{d}t}+m^2$ as in Subsection \ref{Section iteration frame}, we derive the following representation for $V$
\begin{align}\label{repres V(t)}
V(t)= \varepsilon V_{\mathrm{lin}}(t) + \mathrm{e}^{-\alpha_2 t} \int_0^t \mathrm{e}^{(\alpha_2-\alpha_1) s} \int_0^s \mathrm{e}^{\alpha_1 \tau} \, \Gamma(\tau)\left( \int_{\mathbb{R}^n}  |v(\tau,x)|^p \mathrm{d}x\right)^{\beta+1} \mathrm{d}\tau \, \mathrm{d}s,
\end{align} where $\alpha_{1,2}$ are the roots of the quadratic equation $\alpha^2-b\alpha+m^2=0$ and 
\begin{align}\label{def V lin}
V_{\mathrm{lin}}(t) \doteq \begin{cases}\displaystyle{ \frac{\alpha_2\, \mathrm{e}^{-\alpha_1 t}-\alpha_1\,\mathrm{e}^{-\alpha_2 t}}{\alpha_2-\alpha_1} \int_{\mathbb{R}^n} v_0(x) \, \mathrm{d} x +  \,\frac{\mathrm{e}^{-\alpha_1 t}-\mathrm{e}^{-\alpha_2 t}}{\alpha_2-\alpha_1}\int_{\mathbb{R}^n} v_1(x) \, \mathrm{d} x} & \mbox{if} \ b^2>4 m^2,  \\ \displaystyle{\left(1+\tfrac{b}{2}t\right)\mathrm{e}^{-\frac{b}{2} t} \int_{\mathbb{R}^n} v_0(x) \, \mathrm{d} x + t\,  \mathrm{e}^{-\frac{b}{2} t}\int_{\mathbb{R}^n} v_1(x) \, \mathrm{d} x} & \mbox{if} \ b^2=4 m^2.  \end{cases}
\end{align}
From \eqref{repres V(t)} we derive immediately
\begin{align*}
V(t) \gtrsim \mathrm{e}^{-\alpha_2 t} \int_0^t \mathrm{e}^{(\alpha_2-\alpha_1) s} \int_0^s \mathrm{e}^{\alpha_1 \tau} \, \Gamma(\tau) \left(R+cH^{-1}(\mathrm{e}^{H\tau}-1)\right)^{-n(\beta+1)(p-1)}(V(\tau))^{q} \mathrm{d}\tau \, \mathrm{d}s,
\end{align*} where we used H\"older's inequality and the support condition in \eqref{support condition sol anti dS}, and $q=(\beta+1)p$.

Finally, using the inequality $R+\frac cH (\mathrm{e}^{H\tau}-1)\leqslant \left(R+\frac{c}{H}\right)\mathrm{e}^{H\tau}$ for $\tau\geqslant 0$, from the previous inequality we obtain the iteration frame for $V$
\begin{align}\label{Iteration frame V}
V(t)\geqslant C \mathrm{e}^{-\alpha_2 t} \int_0^t \mathrm{e}^{(\alpha_2-\alpha_1) s} \int_0^s \mathrm{e}^{\alpha_1 \tau} \, \Gamma(\tau) \, \mathrm{e}^{-nH(\beta+1)(p-1)\tau}(V(\tau))^{q} \mathrm{d}\tau \, \mathrm{d}s,
\end{align}  where $C=C(n,c,H,\beta,p,R)$ is a suitable positive constant.

Furthermore, \eqref{repres V(t)} provides us a first lower bound estimate for $V$ as well, namely,
\begin{align} \label{1st lb V}
V(t)\geqslant \begin{cases} K_0 \varepsilon \, \mathrm{e}^{-\left(\frac{b}{2}-\frac{1}{2}\sqrt{b^2-4m^2}\right)t} & \mbox{if} \ b^2>4m^2, \\
K_0 \varepsilon \, (1+t) \, \mathrm{e}^{-\frac{b}{2} t} & \mbox{if} \ b^2=4m^2, 
\end{cases} 
\end{align} for any $t\in (0,T)$, where $K_0=K_0(b,m^2,v_0,v_1)$ is a suitable positive and independent of $\varepsilon$ constant. 

Notice that such lower bound for $V$ is completely analogous to the one in \eqref{1st lb U} for $U$, due to the fact that they both follow from the ``linear part'' of $U$ and $V$,  respectively. Nevertheless, for anti-de Sitter spacetime we can exploit the nonlinear term in order to get an alternative lower bound estimate to start the iteration procedure. In the next subsection we will introduce an auxiliary functional that will allow us to derive this further lower bound for $V$.  

\subsection{Lower bound estimate for the nonlinearity}
In the present subsection, we investigate the growth properties of the auxiliary functional given by the following weighted space average of $v$
\begin{align}\label{def V0}
V_0(t)\doteq \int_{\mathbb{R}^n} v(t,x) \Psi(t,x) \, \mathrm{d}x,
\end{align} where the weight function $\Psi=\Psi(t,x; c,H,b,m^2)$ is going to be introduced in few lines and it is chosen as a positive solution of the adjoint homogeneous equation, namely,
\begin{equation} \label{Adjoint Hom Eq}
 \partial_t^2 \Psi -c^2\mathrm{e}^{2Ht} \Delta\Psi-b\partial_t \Psi +m^2\Psi=0. 
\end{equation}
We work with a function $\Psi$ with separable variables, namely, we use the following ansatz $$\Psi(t,x)\doteq \lambda(t;c,H,b,m^2)\Phi(x).$$
As $x$-dependent function we consider the well-known ``eigenfunction'' for the Laplace operator
\begin{align*}
\Phi(x) & \doteq  \mathrm{e}^{x}+\mathrm{e}^{-x}  \qquad  \qquad \ \mbox{if} \ n=1, \\  \Phi(x) & \doteq \int_{\mathbb{S}^{n-1}} \mathrm{e}^{x\cdot \omega} \, \mathrm{d} \sigma_\omega \qquad \mbox{if} \ n\geqslant 2. 
\end{align*} This function has been introduced for the first time in the study of blow-up results for wave models in \cite{YZ06}. The function $\Phi$ is a positive smooth function that satisfies the following crucial properties:
\begin{align}
& \Delta \Phi =\Phi , \label{Laplace Phi =Phi} \\
& \Phi (x) \sim c_n |x|^{-\frac{n-1}{2}}  \mathrm{e}^x \qquad \mbox{as} \ |x|\to \infty, \label{Asymptotic behavior Phi}
\end{align} for some suitable positive constant $c_n$.
Also, in order to get a solution of \eqref{Adjoint Hom Eq}, we have to determine $\lambda=\lambda(t;c,H,b,m^2)$ such that
\begin{equation} \label{Lambda Eq}
 \frac{\mathrm{d}^2 \lambda}{\mathrm{d}t^2} -b \, \frac{\mathrm{d} \lambda}{\mathrm{d}t}  +\left(m^2- c^2\mathrm{e}^{2Ht}\right)\lambda =0.
\end{equation} For the sake of readability, in what follows we skip the dependence of $\lambda$ on $(c,H,b,m^2)$  in the notations. 
 Let us perform the change of variables $\tau=\frac cH \mathrm{e}^{Ht}$. Then,
\begin{align*}
\frac{\mathrm{d}^2 \lambda}{\mathrm{d}t^2}=  H^2\tau^2 \frac{\mathrm{d}^2 \lambda}{\mathrm{d}\tau^2}+H^2\tau \frac{\mathrm{d} \lambda}{\mathrm{d}\tau} \qquad   \mbox{and} \qquad \frac{\mathrm{d} \lambda}{\mathrm{d}t} =  H\tau \frac{\mathrm{d} \lambda}{\mathrm{d}\tau} .
\end{align*} So, with respect to $\tau$ the function $\lambda$ satisfies the equation
\begin{equation} \label{Lambda Eq tau}
\tau^2 \frac{\mathrm{d}^2 \lambda}{\mathrm{d}\tau^2} +\left(-\frac bH+1\right)\tau \, \frac{\mathrm{d} \lambda}{\mathrm{d}\tau}  +\left(\frac{m^2}{H^2}- \tau^2\right)\lambda =0.
\end{equation} Next, we carry out the transformation $\lambda(\tau)=\tau^\rho \eta(\tau)$, with $\rho$ real parameter to be determined. 
By straightforward computations, we have that $\lambda$ solves \eqref{Lambda Eq tau} if and only if 
\begin{equation} \label{Eta Eq}
 \tau^2 \frac{\mathrm{d}^2 \eta}{\mathrm{d}\tau^2} +\left(2\rho-\frac{b}{H}+1\right) \tau \, \frac{\mathrm{d} \eta}{\mathrm{d}\tau}  +\left(\rho\left(\rho-\frac{b}{H}\right)+\frac{m^2}{H^2}- \tau^2\right)\eta =0.
\end{equation}
If we choose $\rho\doteq \frac{b}{2H}$, then, \eqref{Eta Eq} can be rewritten as the following modified Bessel equation
\begin{equation}\label{Bessel Eq}
\tau^2 \frac{\mathrm{d}^2 \eta}{\mathrm{d}\tau^2} + \tau \, \frac{\mathrm{d} \eta}{\mathrm{d}\tau}  -\left[\frac{1}{4H^2}\left(b^2-4m^2\right)+ \tau^2\right]\eta =0.
\end{equation} Setting $\nu \doteq \frac{1}{2H} (b^2-4m^2)^{1/2}$, a complete system of independent solutions to \eqref{Bessel Eq} is given by $\mathrm{I}_\nu$ and $\mathrm{K}_\nu$ (modified Bessel functions of the first and second kind, respectively, of order $\nu$). For further details on the properties of $\mathrm{I}_\nu$ and $\mathrm{K}_\nu$ that will be used in this subsection, we address the reader to \cite[Chapter 10]{OLBC10}. Let us recall the asymptotic behavior of $\mathrm{I}_\nu$ and $\mathrm{K}_\nu$ 
 for large values
\begin{align*}
& \mathrm{I}_\nu(\tau) \sim (2\pi\tau)^{-1/2} \mathrm{e}^\tau \quad \mbox{and} \quad \mathrm{K}_\nu(\tau) \sim \sqrt{\frac{\pi}{2}} \, \tau^{-1/2} \mathrm{e}^{-\tau} \qquad \mbox{as} \ \tau \to \infty.
\end{align*} Moreover, for $\nu\geqslant 0$ the function $\mathrm{I}_\nu(\tau)$ has no real zero excluding $\tau=0$ when $\nu>0$, and, similarly, for $\nu\geqslant 0$ the function $\mathrm{K}_\nu(\tau)$ has no real zero.

Hereafter, we set (neglecting the unessential multiplicative constant)
\begin{align}\label{def lambda anti dS}
\lambda(t;c,H,b,m^2)\doteq \mathrm{e}^{\frac{b}{2}t} \, \mathrm{K}_\nu\left(\tfrac cH \mathrm{e}^{Ht}\right).
\end{align} By using the previous recalled asymptotic behavior of $\mathrm{K}_\nu$ for large arguments and the fact that $\mathrm{K}_\nu$ has no real zero, we may consider the following uniform estimate
\begin{align}\label{Asymp Lambda}
\lambda_0 \, \mathrm{e}^{\frac{1}{2}(b-H)t}\exp\left(-\frac{c}{H}\mathrm{e}^{Ht}\right) \leqslant \lambda(t)\leqslant \Lambda_0 \,  \mathrm{e}^{\frac{1}{2}(b-H)t}\exp\left(-\frac{c}{H}\mathrm{e}^{Ht}\right) \qquad \mbox{for any} \ t\geqslant 0,
\end{align} for some positive constants $\lambda_0=\lambda_0(c,H,b,m^2),\Lambda_0=\Lambda_0(c,H,b,m^2)$.

By using the uniform estimate \eqref{Asymp Lambda}, we can now derive a lower bound for the functional $V_0$. As $\lambda$ and $\Phi$ are nonnegative functions, also $\Psi$ is nonnegative. Therefore, plugging $\Psi$ in \eqref{def weak sol int rel anti dS}, we get
\begin{align*}
0 & \leqslant  \int_0^t \Gamma(s) \left(\int_{\mathbb{R}^n}  |v(s,y)|^p  \, \mathrm{d}y\right)^\beta \int_{\mathbb{R}^n}  |v(s,x)|^p  \Psi(s,x) \, \mathrm{d}x \, \mathrm{d}s  \\
& = \int_0^t \int_{\mathbb{R}^n} v(s,x)\left(  \Psi_{ss}(s,x)- c^2\mathrm{e}^{2Hs} \Delta  \Psi(s,x)- b  \Psi_s(s,x) +m^2  \Psi(s,x) \right) \, \mathrm{d}x \, \mathrm{d}s \\
 & \qquad + \int_{\mathbb{R}^n} \left(\partial_t v(s,x)\Psi(s,x)-v(s,x) \Psi_s(s,x)+b v(s,x)\Psi(s,x)\right) \mathrm{d}x \, \Big|^{s=t}_{s=0} \\
 & = V'_0(s)+bV_0(s)-2 \frac{\lambda'(s)}{\lambda(s)}V_0(s) \Big|^{s=t}_{s=0},
\end{align*} where in the last step we used \eqref{Adjoint Hom Eq} and the trivial relations
\begin{align*}
 V'_0(s) & = \int_{\mathbb{R}^n} \left(\partial_t u(s,x)\Psi(s,x)+u(s,x) \Psi_s(s,x)\right) \mathrm{d}x   , 
  & \frac{\lambda'(s)}{\lambda(s)}V_0(s) = \int_{\mathbb{R}^n} u(s,x)\Psi_s(s,x)\, \mathrm{d}x   .
\end{align*} Let us remark that 
\begin{align*}
\frac{\mathrm{d}}{\mathrm{d}s}\left(\frac{\mathrm{e}^{bs}}{\lambda^2(s)} V_0(s)\right)= \frac{\mathrm{e}^{bs}}{\lambda^2(s)}\left[V'_0(s)+\left(b-2 \frac{\lambda'(s)}{\lambda(s)}\right)V_0(s) \right].
\end{align*} Hence, the previous inequality implies
\begin{align}
\lambda^2(t)\, \mathrm{e}^{-bt}\frac{\mathrm{d}}{\mathrm{d}t}\left(\frac{\mathrm{e}^{bt}}{\lambda^2(t)} V_0(t)\right) & \geqslant V'_0(0)+bV_0(0)-2 \frac{\lambda'(0)}{\lambda(0)}V_0(0) \notag\\
& = \varepsilon \underbrace{ \int_{\mathbb{R}^n} \left[\lambda(0) v_1(x)+(b\lambda(0)-\lambda'(0))v_0(x) \right] \Phi(x) \,\mathrm{d}x}_{\doteq \mathcal{I}[v_0,v_1]}. \label{ODI V0}
\end{align} 
Using the recursive relation 
\begin{align*}
\frac{\partial \mathrm{K}_\nu}{\partial z}(z) = -  \mathrm{K}_{\nu+1}(z) +\frac{\nu}{z}  \mathrm{K}_{\nu }(z)
\end{align*} (cf. \cite[Section 10.29]{OLBC10}), it follows
\begin{align*}
\lambda'(t)& =\frac{b}{2}\lambda(t)+c \, \mathrm{e}^{(\frac{b}{2}+H)t} \mathrm{K}'_\nu\left(\tfrac cH \mathrm{e}^{Ht}\right) \\
& =\frac{b}{2}\lambda(t)+c\, \mathrm{e}^{(\frac{b}{2}+H)t} \left[  -\mathrm{K}_{\nu+1}\left(\tfrac cH \mathrm{e}^{Ht}\right) +\frac{\nu H}{c}\mathrm{e}^{-Ht} \mathrm{K}_\nu\left(\tfrac cH \mathrm{e}^{Ht}\right)\right] \\
& =\left(\frac{b}{2}+\frac{1}{2}\sqrt{b^2-4m^2}\right)\lambda(t)-c\, \mathrm{e}^{(\frac{b}{2}+H)t} \mathrm{K}_{\nu+1}\left(\tfrac cH \mathrm{e}^{Ht}\right) 
\end{align*} which implies in turn
\begin{align*}
b\lambda(t) - \lambda'(t) & = \left(\frac{b}{2}-\frac{1}{2}\sqrt{b^2-4m^2}\right)\lambda(t)+c\, \mathrm{e}^{(\frac{b}{2}+H)t} \mathrm{K}_{\nu+1}\left(\tfrac cH \mathrm{e}^{Ht}\right) \geqslant 0.
\end{align*} 
 Therefore, assuming $v_0,v_1$ nonnegative and nontrivial, in particular, we have that $\mathcal{I}[v_0,v_1]>0$.
From \eqref{ODI V0}, we find
\begin{align*}
 & \frac{\mathrm{d}}{\mathrm{d}t}\left(\frac{\mathrm{e}^{bt}}{\lambda^2(t)} V_0(t)\right)  \geqslant  \varepsilon \mathcal{I}[v_0,v_1] \frac{\mathrm{e}^{bt}}{\lambda^2(t)}, 
 \end{align*}
 from which it follows
 \begin{align} V_0(t)\geqslant  \frac{V_0(0)}{\lambda^2(0)} \mathrm{e}^{-bt} \lambda^2(t)  +  \varepsilon \mathcal{I}[v_0,v_1]  \, \mathrm{e}^{-bt}\lambda^2(t)  \int_0^t \frac{\mathrm{e}^{bs}}{\lambda^2(s)} \mathrm{d}s. \label{lower bound F0 general}
\end{align} By \eqref{Asymp Lambda}, it results
\begin{align*}
\frac{V_0(0)}{\lambda^2(0)} \mathrm{e}^{-bt} \lambda^2(t) \geqslant \varepsilon \frac{ \lambda_0^2}{\lambda(0)} \left( \int_{\mathbb{R}^n}v_0(x)  \Phi(x) \, \mathrm{d}x \right) \mathrm{e}^{-H t}\exp\left(-\frac{2c}{H}\mathrm{e}^{H t}\right)
\end{align*} and 
\begin{align*}
\mathrm{e}^{-bt}\lambda^2(t)  \int_0^t \frac{\mathrm{e}^{bs}}{\lambda^2(s)} \, \mathrm{d}s  & \geqslant \frac{\lambda_0^2}{\Lambda_0^2} \,  \mathrm{e}^{-H t}\exp\left(-\frac{2c}{H}\mathrm{e}^{H t}\right) \int_0^t \mathrm{e}^{H s}\exp\left(\frac{2c}{H}\mathrm{e}^{H s}\right) \mathrm{d}s \\
& = \frac{\lambda_0^2}{2c \Lambda_0^2} \,  \mathrm{e}^{-H t}\exp\left(-\frac{2c}{H}\mathrm{e}^{H t}\right) \left[\exp\left(\frac{2c}{H}\mathrm{e}^{H s}\right) \right]_{s=0}^{s=t} \\
& = \frac{\lambda_0^2}{2c \Lambda_0^2} \,  \mathrm{e}^{-H t}\left(1-\exp\left(-\frac{2c}{H}\mathrm{e}^{H t}\right) \right).
\end{align*} Using the first of the two previous estimates for $t$ in a neighborhood of $0$ and the second one for $t$ away from zero, we derive the following lower bound estimate
\begin{align}\label{1st lb V0}
V_0(t)\gtrsim \varepsilon \, \mathrm{e}^{-H t} \qquad \mbox{for} \ t\in [0,T).
\end{align}

The next step is to derive a lower bound estimate for $\|v(t,\cdot)\|_{L^p(\mathbb{R}^n)}^p$ by using the previous lower bound for $V_0$. By H\"older's inequality and using the support condition \eqref{support condition sol anti dS}, for $t\geqslant 0$ it follows
\begin{align*}
V_0(t) \leqslant  \|v(t,\cdot)\|_{L^p(\mathbb{R}^n)} \|\Psi(t,\cdot)\|_{L^{p'}\left(B_{R+A_{\mathrm{AdS}}(t)}\right)}.
\end{align*}  Consequently, if we get an upper bound for $\|\Psi(t,\cdot)\|_{L^{p'}\left(B_{R+A_{\mathrm{AdS}}(t)}\right)}$, then, we may obtain the desired estimate by using
\begin{align}\label{intermediate lower bound |v|^p}
\int_{\mathbb{R}^n} |v(t,x)|^p \, \mathrm{d} x \geqslant (V_0(t))^p \, \|\Psi(t,\cdot)\|_{L^{p'}\left(B_{R+A_{\mathrm{AdS}}(t)}\right)}^{-p}.
\end{align}
 Repeating the same computations as in \cite[Section 3]{PalRei19}, we get
\begin{align*}
\|\Psi(t,\cdot)\|_{L^{p'}\left(B_{R+A_{\mathrm{AdS}}(t)}\right)} & = \left(\int_{B_{R+A_{\mathrm{AdS}}(t)}} (\Psi(t,x))^{p'} \mathrm{d}x\right)^{1/p'} \leqslant \lambda(t) \! \left(\int_{B_{R+A_{\mathrm{AdS}}(t)}} (\Phi(x))^{p'} \mathrm{d}x\right)^{1/p'}  \\ & \lesssim \lambda(t) \exp (R+A_{\mathrm{AdS}}(t)) (R+A_{\mathrm{AdS}}(t))^{(n-1)\big(\frac{1}{2}-\frac{1}{p}\big)}.
\end{align*} Thus, from \eqref{Asymp Lambda} for $t\geqslant 0$ we get 
\begin{align}
\|\Psi(t,\cdot)\|_{L^{p'}\left(B_{R+A_{\mathrm{AdS}}(t)}\right)} & \lesssim  \Lambda_0 \,  \mathrm{e}^{\frac{1}{2}(b-H)t}\exp\!\left(-\tfrac{c}{H}\mathrm{e}^{Ht}\right)  \exp (R+A_{\mathrm{AdS}}(t)) (R+A_{\mathrm{AdS}}(t))^{(n-1)\big(\frac{1}{2}-\frac{1}{p}\big)} \notag \\
& \lesssim    \mathrm{e}^{\frac{1}{2}(b-H)t} \left(R+\tfrac cH \left(\mathrm{e}^{Ht} -1\right)\right)^{(n-1)(\frac{1}{2}-\frac{1}{p})}. \label{upper bound integral Psi anti dS}
\end{align}  Combining \eqref{1st lb V0}, \eqref{intermediate lower bound |v|^p} and \eqref{upper bound integral Psi anti dS}, we conclude
\begin{equation} \label{Lower bound |v|^p int}
\int_{\mathbb{R}^n} |v(t,x)|^p \, \mathrm{d} x  \gtrsim \varepsilon^p \, \mathrm{e}^{-\frac{1}{2}\left(b+H\right)p t} \left(R+\tfrac cH\left(\mathrm{e}^{Ht} -1\right)\right)^{(n-1)\big(1-\frac{p}{2}\big)},
\end{equation} for $t\geqslant 0$.

Finally, plugging \eqref{Lower bound |v|^p int} in \eqref{repres V(t)}, for $t\geqslant t_0 >0$ we get 
\begin{align}\label{1st lb V from V0}
V(t)\geqslant K_1  \varepsilon^{(\beta+1)p} \, \mathrm{e}^{\big[\varrho-(b+H)(\beta+1)\frac{p}{2} +(n-1)H(\beta+1)\big(1-\frac{p}{2}\big)\big] t} t^{\varsigma},
\end{align} where the constant $K_1=K_1(c,H,b,m^2,v_0,v_1,p,\beta,\mu,\varrho,\varsigma,t_0)>0$ is independent of $\varepsilon$ and $t_0$ will be fixed time by time depending on the slicing procedure that we will apply in each case.

\begin{remark} In \eqref{def lambda anti dS} we considered a modified Bessel function of the second kind of order $\nu$. Nonetheless, if we had chosen a modified Bessel function of the first kind of order $\nu$ instead, namely, defining $\lambda(t)\doteq \mathrm{e}^{\frac{b}{2}t} \, \mathrm{I}_\nu\left(\tfrac cH \mathrm{e}^{Ht}\right)$ in place of \eqref{def lambda anti dS}, the final outcome (meaning the lower bound bound for the space integral of the power nonlinearity) would have been the same. More in detail, the lower bound for $V_0$ in \eqref{1st lb V0} would have been $V_0(t)\gtrsim \varepsilon\,  \mathrm{e}^{-tH}\exp\left(2cH^{-1} \mathrm{e}^{Ht}\right)$. In particular, this better lower bound for $V_0$ would have been a consequence of the asymptotic behavior for large argument of the factor involving $\mathrm{I}_\nu$ in this alternative definition of $\lambda$. Moreover, this lower bound would have followed from the estimate of the first term in \eqref{lower bound F0 general}, since the contribute from the integral term in \eqref{lower bound F0 general} would have been weaker in this case. Nevertheless, the better behavior of $\lambda$ would have influenced the estimate \eqref{upper bound integral Psi anti dS} as well, providing eventually exactly the same estimate for  $\|v(t,\cdot)\|_{L^p(\mathbb{R}^n)}^p$ as the one in \eqref{Lower bound |v|^p int}.
\end{remark}

\begin{remark} \label{Remark V0 for de Sitter} Clearly even in the case of de Sitter spacetime we could have considered a weighted functional analogous to the one in \eqref{def V0}. However, the lower bound estimates corresponding to the one in \eqref{1st lb V from V0} in this case would have been exactly the same as we have plugged \eqref{1st lb U} in \eqref{Iteration frame}. In other words, we would have start the sequence of lower bound estimates in \eqref{lb estimate U step j exponential growth} from $j=1$ rather than from $j=0$, but obviously the final outcome would have be unaltered.
\end{remark}

\subsection{Comparison of the first lower bound}
\label{Subsection comparison 1st lb estimates}

In order to prove Theorems  \ref{Thm anti dS lin lb} - \ref{Thm anti dS nlin lb} 
it is crucial to understand which first lower bound estimate for $V$ between \eqref{1st lb V} and \eqref{1st lb V from V0} has the dominant role. 

We remark that \eqref{1st lb V} and \eqref{1st lb V from V0} cannot be directly compared in their current forms due to the different orders for $\varepsilon$ in their right-hand sides. Then, we plug \eqref{1st lb V} in \eqref{Iteration frame V} obtaining 
\begin{align} \label{2nd lb V} 
V(t)\gtrsim \varepsilon^{(\beta+1)p} (t-L_2)^{\varsigma_+ +q b_0} t^{-\varsigma_-} \mathrm{e}^{\left(\varrho -nH(\beta+1)(p-1)-\frac{1}{2}\left(b-\sqrt{b^2-4m^2}\right)(\beta+1)p\right)t}
\end{align} for $t\geqslant L_2$ (see next section for the derivation of this inequality and the definition of $L_2$).
Now, we can compare the multiplicative coefficient in the exponential term in \eqref{1st lb V from V0} and in \eqref{2nd lb V}. In particular, we have that the coefficient in the exponential term in \eqref{1st lb V from V0} is dominant over the one in \eqref{2nd lb V} provided that 
\begin{align*}
\varrho -nH(\beta+1)(p-1)-\tfrac{1}{2}\left(b-\sqrt{b^2-4m^2}\right)q< \varrho-(b+H)\tfrac{q}{2} +(n-1)H(\beta+1)\left(1-\tfrac{p}{2}\right)
\end{align*} and, by straightforward computations, we have that  the previous inequality is equivalent to \eqref{1st lb for V from V0 is dominant}. 
 Notice that \eqref{1st lb for V from V0 is dominant} is exactly the condition that provides a $\varrho_{\mathrm{crit}}(n,H,b,m^2,\beta,p)$ given by \eqref{def rho crit V0}. Hence, when \eqref{1st lb for V from V0 is dominant} holds we shall use \eqref{1st lb V from V0} as first lower bound estimate for $V$, while when \eqref{1st lb for V from lin part is dominant} holds, we shall use \eqref{1st lb V}. Notice that even in the limit case $\frac{n}{2}-\frac{\sqrt{b^2-4m^2}}{2H}=\frac{1}{p}$ we use \eqref{1st lb V}, since in the case $b^2=4m^2$ a slight improvement of polynomial type is included, whilst this does not happen in \eqref{1st lb V from V0}. So far we discussed which among \eqref{1st lb V} and \eqref{1st lb V from V0} is better to start the iteration argument depending on the values of $n,H,b,m^2$ and $p$. Clearly, the reasons behind the actual definition of $\varrho_{\mathrm{crit}}$ either as in \eqref{def rho crit lin part} or as in \eqref{def rho crit V0} will be clarified in the proofs of Theorem \ref{Thm anti dS lin lb} and Theorem \ref{Thm anti dS nlin lb}, respectively, by the corresponding iteration arguments.

We point out that \eqref{1st lb for V from V0 is dominant} is never satisfied for $n\leqslant  N_0$, where 
\begin{align}\label{def N0 dim}
N_0= N_0(H,b,m^2) \doteq \frac{\sqrt{b^2-4m^2}}{H}
\end{align} since the left-hand side of \eqref{1st lb for V from V0 is dominant} is negative for $n$ in this range, while for $n\geqslant  2+N_0$ the condition in  \eqref{1st lb for V from V0 is dominant} is always fulfilled since since the left-hand side is greater than or equal to $1$. Finally, for $n\in \left(N_0, 2+N_0\right)$ we have that \eqref{1st lb for V from V0 is dominant} is true if and only if 
\begin{align*}
p>\widetilde{p}(n,H,b,m^2)\doteq \frac{2H}{nH -\sqrt{b^2-4m^2}}.
\end{align*}
Summarizing, for $n\leqslant N_0$ or $n\in (N_0, 2+N_0)$ and $1<p\leqslant \widetilde{p}$ we will use \eqref{1st lb V} to star the iteration argument, while for  $n\in (N_0, 2+N_0)$ and $p>\widetilde{p}$ or $n\geqslant 2+N_0$ we will employ \eqref{1st lb V from V0} as staring point for the iteration procedure. 

Finally, we emphasize that the previous conditions for the employment of either \eqref{1st lb V} or \eqref{1st lb V from V0} (which correspond to \eqref{1st lb for V from lin part is dominant} and \eqref{1st lb for V from V0 is dominant}, respectively) are completely independent of $\beta$.

\subsection{Proof of Theorems \ref{Thm anti dS lin lb}, \ref{Thm anti dS lin lb poly growth} and \ref{Thm anti dS lin lb log growth}}

As we explained in Subsection \ref{Subsection comparison 1st lb estimates}, for $n\leqslant N_0$ or $n\in (N_0, 2+N_0)$ and $1<p\leqslant \widetilde{p}$  (that is, when \eqref{1st lb for V from lin part is dominant} holds) we employ \eqref{1st lb V} as first lower bound estimate for $V$. By plugging the explicit expression for the factor $\Gamma$ into \eqref{Iteration frame V}, we find
\begin{align}\label{Iteration frame V lin}
V(t)\geqslant C \mu (1+t)^{-\varsigma_-}  \mathrm{e}^{-\alpha_2 t} \int_0^t \mathrm{e}^{(\alpha_2-\alpha_1) s} \int_0^s \mathrm{e}^{\left[\alpha_1+\varrho-  nH(\beta+1)(p-1)\right]\tau} \, \tau^{\varsigma_+}\, (V(\tau))^{q} \mathrm{d}\tau \, \mathrm{d}s
\end{align} for $t\geqslant 0$, where $\varsigma_+,\varsigma_-$ denote the positive and the negative part of $\varsigma$, respectively.

We recognize that the iteration frame in \eqref{Iteration frame V lin} is formally identical to the one in \eqref{Iteration frame exponential growth} for the functional $U$. Moreover, the lower bound estimates for $U$ and $V$ in \eqref{1st lb U} and \eqref{1st lb V} are completely analogous (up to the multiplicative constants that depend on the Cauchy data) and $$\varrho_{\mathrm{crit}}(n,H,b,m^2,\beta,p)-nH(\beta+1)(p-1)=r_{\mathrm{crit}}(b,m^2,\beta,p).$$ Hence, the proof of Theorem \ref{Thm anti dS lin lb} is completely similar to the one of Theorem \ref{Thm exponential growth}, provided that we consider as first parameter that characterizes the slicing procedure $$\ell_0\doteq \max\left\{(\alpha_1+\varrho -nH(\beta+1)(p-1))^{-1},(\alpha_2+\varrho -nH(\beta+1)(p-1))^{-1}\right\}$$ and then $$L_j\doteq \ell_0 \prod_{k=1}^j \left(1+q^{-\frac{k}{2}}\right) \qquad \mbox{for any} \ j\in\mathbb{N}^*.$$

On the other hand, Theorems \ref{Thm anti dS lin lb poly growth} and \ref{Thm anti dS lin lb log growth} can be proved analogously as Theorems \ref{Thm polynomial growth} and \ref{Thm logarithmic growth} by working with the functional $\mathcal{V}(t)\doteq \mathrm{e}^{\alpha_1 t}V(t)$. Indeed, setting $\alpha_1=\frac{1}{2}(b-\sqrt{b^2-4m^2})$ and proceeding as in Subsection \ref{Subsection polynomial growth}, from \eqref{Iteration frame V} we get 
\begin{align}\label{Iteration frame mathcal V}
\mathcal{V}(t)\geqslant C \mu \mathrm{e}^{(\alpha_1-\alpha_2)t} \int_0^t \mathrm{e}^{(\alpha_2-\alpha_1)s}  \int_0^s (1+\tau)^\varsigma (\mathcal{V}(\tau))^{q} \mathrm{d}\tau \, \mathrm{d}s,
\end{align} for $t\geqslant 0$. So, $\mathcal{V}$ has an iteration frame formally identical to the one for $\mathcal{U}$ in \eqref{Iteration frame mathcal U}. Besides, $\mathcal{V}$ satisfies a completely analogous lower bound estimate as the one in \eqref{1st lb mathcal U} for $\mathcal{U}$. By following the same approaches as for the proofs in Subsections  \ref{Subsection polynomial growth} and \ref{Subsection log growth} we conclude the validity of Theorems \ref{Thm anti dS lin lb poly growth} and \ref{Thm anti dS lin lb log growth}. 

\subsection{Proof of Theorem \ref{Thm anti dS nlin lb}}

According to Subsection \ref{Subsection comparison 1st lb estimates}, for $n\in (N_0, 2+N_0)$ and $p> \widetilde{p}$ or $n\geqslant 2+N_0$  (that is, when \eqref{1st lb for V from V0 is dominant} holds) it is appropriate to consider \eqref{1st lb V from V0} as first lower bound estimate for $V$.

Our first goal is to prove the sequence of lower bound estimates 
\begin{align} \label{lb estimate V step j exponential growth}
V(t)\geqslant C_j \mathrm{e}^{a_j t}\mathrm{e}^{-\gamma_j t} (t-L_{2j})^{b_j}(1+t)^{-\beta_j} \quad \mbox{for} \ t\geqslant L_{2j} \ \mbox{and for any} \ j\in\mathbb{N},
\end{align} where $\{C_j\}_{j\in\mathbb{N}}$, $\{a_j\}_{j\in\mathbb{N}}$, $\{\gamma_j\}_{j\in\mathbb{N}}$, $\{b_j\}_{j\in\mathbb{N}}$, $\{\beta_j\}_{j\in\mathbb{N}}$ are sequences of nonnegative real numbers to be determined iteratively and the parameters characterizing the slicing procedure $\{L_j\}_{j\in\mathbb{N}}$ are given formally by \eqref{def Lj exponential growth}, however,  with the following modifications in the definition of the parameters $\{\ell_k\}_{k\in\mathbb{N}}$ 
\begin{align*}
\ell_0 &\doteq \max\left\{(A_0+\alpha_1)^{-1},(A_0+\alpha_2)^{-1}\right\}, \qquad
\ell_k \doteq 1+ ((\beta+1)p)^{-\frac{k}{2}} \quad \mbox{for} \ k\geqslant 1,
\end{align*} where 
\begin{align*}
A_0 
\doteq \varrho-\varrho_{\mathrm{crit}}(n,H,b,m^2,\beta,p)+\tfrac{1}{2}(b+nH)(q-1)+H (\beta+1).
\end{align*} Moreover, setting
\begin{align*}
A_1 
\doteq nH(q-1)+\tfrac{H}{p},
\end{align*} we can represent the part of the coefficient for the exponential term in the $\tau$-integral in \eqref{Iteration frame V lin} given by
\begin{align}
\varrho -nH(\beta+1)(p-1) & = \varrho-\varrho_{\mathrm{crit}}+\tfrac{1}{2}(b+nH)(q-1)+H (\beta+1)-nH(q-1)-\tfrac{H}{p} \notag \\ &= A_0-A_1\label{A0-A1}
\end{align} as difference between the two positive quantities $A_0,A_1$. As we will see in the iteration frame, the splitting \eqref{A0-A1} is the reason for the previous choice of $\ell_0$.

Clearly, \eqref{1st lb V from V0} implies the validity of \eqref{lb estimate V step j exponential growth} for $j=0$ provided that $$C_0\doteq K_1 \varepsilon^q, \  a_0\doteq \varrho- \varrho_{\mathrm{crit}}+\tfrac{nH}{2}, \ \gamma_0\doteq \tfrac b2 +\tfrac Hp, \ b_0\doteq \varsigma_+, \ \beta_0\doteq \varsigma_-,$$
where we applied the decomposition $ a_0-\gamma_0=\varrho -(b+H)\tfrac q2 +(n-1)H(\beta+1)\left(1-\tfrac p2\right)$ for the coefficient of the exponential term in \eqref{1st lb V from V0}.

Plugging \eqref{Gamma factor exponential growth anti dS nlin lb} in \eqref{Iteration frame V} and using \eqref{A0-A1}, we get the following iteration frame
\begin{align}
V(t)& \geqslant C \mu (1+t)^{-\varsigma_-}  \mathrm{e}^{-\alpha_2 t} \int_0^t \mathrm{e}^{(\alpha_2-\alpha_1) s} \int_0^s \mathrm{e}^{\left[\alpha_1+\varrho-  nH(\beta+1)(p-1)\right]\tau} \, \tau^{\varsigma_+}\, (V(\tau))^{q} \mathrm{d}\tau \, \mathrm{d}s \notag \\
& \geqslant C \mu (1+t)^{-\varsigma_-}  \mathrm{e}^{-(\alpha_2+A_1) t} \int_0^t \mathrm{e}^{(\alpha_2-\alpha_1) s} \int_0^s \mathrm{e}^{(\alpha_1+A_0)\tau} \, \tau^{\varsigma_+}\, (V(\tau))^{q} \mathrm{d}\tau \, \mathrm{d}s
\label{Iteration frame V nlin lb}
\end{align} for $t\geqslant 0$.

Let us prove the induction step: assuming \eqref{lb estimate V step j exponential growth} satisfied for some $j\geqslant 0$, we prove it for $j+1$. Plugging \eqref{lb estimate V step j exponential growth} in \eqref{Iteration frame V nlin lb}, we get
\begin{align} \label{lb V step j+i intermediate}
V(t)& \geqslant \frac{\mu \, C C_j^q}{(1+t)^{\varsigma_-+q\beta_j}} \,  \mathrm{e}^{-(\alpha_2+A_1+q\gamma_j) t} \int_{L_{2j}}^t \mathrm{e}^{(\alpha_2-\alpha_1) s} \int_{L_{2j}}^s \mathrm{e}^{(\alpha_1+A_0+a_jq)\tau}  (\tau-L_{2j})^{\varsigma_++qb_j}\,  \mathrm{d}\tau \, \mathrm{d}s
\end{align} for $t\geqslant L_{2j}$. Let us apply a slicing procedure to the $\tau$-integral. For $t\geqslant L_{2j+1}$ we may estimate
\begin{align*}
 & \int_{L_{2j}}^t \mathrm{e}^{(\alpha_2-\alpha_1) s} \int_{L_{2j}}^s \mathrm{e}^{(\alpha_1+A_0+a_jq)\tau}  (\tau-L_{2j})^{\varsigma_++qb_j}\,  \mathrm{d}\tau \, \mathrm{d}s \\ & \quad  \geqslant  \int_{L_{2j+1}}^t \mathrm{e}^{(\alpha_2-\alpha_1) s} \int_{\tfrac{L_{2j}s}{L_{2j+1}}}^s \mathrm{e}^{(\alpha_1+A_0+a_jq)\tau}  (\tau-L_{2j})^{\varsigma_++qb_j}\,  \mathrm{d}\tau \, \mathrm{d}s 
 \\ & \quad  \geqslant  \ell_{2j+1}^{-(\varsigma_++qb_j)}  \int_{L_{2j+1}}^t \mathrm{e}^{(\alpha_2-\alpha_1) s}  (s-L_{2j+1})^{\varsigma_++qb_j} \int_{\tfrac{L_{2j}s}{L_{2j+1}}}^s \mathrm{e}^{(\alpha_1+A_0+a_jq)\tau} \,  \mathrm{d}\tau \, \mathrm{d}s \\
 & \quad  =  \frac{\ell_{2j+1}^{-(\varsigma_++qb_j)}}{ A_0+\alpha_1+q a_j} \int_{L_{2j+1}}^t  (s-L_{2j+1})^{\varsigma_++qb_j} \mathrm{e}^{(\alpha_2+A_0+a_jq)s} \left(1-\mathrm{e}^{-(\alpha_1+A_0+a_jq)(1-1/\ell_{2j+1})s}\right) \, \mathrm{d}s.
\end{align*} For $s\geqslant L_{2j+1}$ we can estimate
\begin{align*}
1-\mathrm{e}^{-(\alpha_1+A_0+a_jq)(1-1/\ell_{2j+1})s} & \geqslant 1-\mathrm{e}^{-(\alpha_1+A_0+a_jq)(\ell_{2j+1}-1)L_{2j}}   \geqslant 1-\mathrm{e}^{-(\alpha_1+A_0+a_jq)(\ell_{2j+1}-1)\ell_0}  \\
& \geqslant 1-\mathrm{e}^{-(\alpha_1+A_0)(\ell_{2j+1}-1)\ell_0}   \geqslant 1-\mathrm{e}^{-(\ell_{2j+1}-1)} \geqslant q^{-(2j+1)}(q-\tfrac 12),
\end{align*} where in the last step we used \eqref{remainder exponential term}, therefore, for $t\geqslant L_{2j+1}$ we have
\begin{align*}
& \int_{L_{2j}}^t \mathrm{e}^{(\alpha_2-\alpha_1) s} \int_{L_{2j}}^s \mathrm{e}^{(\alpha_1+A_0+a_jq)\tau}  (\tau-L_{2j})^{\varsigma_++qb_j}\,  \mathrm{d}\tau \, \mathrm{d}s \\
 & \qquad  \geqslant \frac{(q-\tfrac 12) \, \ell_{2j+1}^{-(\varsigma_++qb_j)} q^{-(2j+1)}}{ A_0+\alpha_1+q a_j} \int_{L_{2j+1}}^t  (s-L_{2j+1})^{\varsigma_++qb_j} \mathrm{e}^{(\alpha_2+A_0+a_jq)s} \, \mathrm{d}s.
\end{align*} Repeating a similar estimate for the $s$-integral, for $t\geqslant L_{2j+2}$ we obtain
\begin{align*}
& \int_{L_{2j}}^t \mathrm{e}^{(\alpha_2-\alpha_1) s} \int_{L_{2j}}^s \mathrm{e}^{(\alpha_1+A_0+a_jq)\tau}  (\tau-L_{2j})^{\varsigma_++qb_j}\,  \mathrm{d}\tau \, \mathrm{d}s \\
 & \qquad  \geqslant \frac{(q-\tfrac 12)^2 \big(\ell_{2j+1} \ell_{2j+2}\big)^{-(\varsigma_++qb_j)} q^{-(4j+3)}}{ (A_0+\alpha_1+q a_j)(A_0+\alpha_2+q a_j)} (t-L_{2j+2})^{\varsigma_++qb_j} \mathrm{e}^{(\alpha_2+A_0+a_jq)t}. 
\end{align*}
Combining the previous inequality and \eqref{lb V step j+i intermediate}, for $t\geqslant L_{2j+2}$ we find
\begin{align*}
V(t)&\geqslant     \frac{\mu \, C C_j^q (q-\tfrac 12)^2 \big(\ell_{2j+1} \ell_{2j+2}\big)^{-(\varsigma_++qb_j)} }{ (A_0+\alpha_1+q a_j)(A_0+\alpha_2+q a_j) q^{4j+3}} \, \mathrm{e}^{[(A_0+q a_j)-(A_1+q\gamma_j) ]t}  \frac{(t-L_{2j+2})^{\varsigma_++qb_j}}{(1+t)^{\varsigma_-+q\beta_j} }
\end{align*} which is \eqref{lb estimate V step j exponential growth} for $j+1$ provided that
\begin{align*}
C_{j+1}& \doteq  \frac{\mu \, C  (q-\tfrac 12)^2 C_j^q }{ (A_0+\alpha_1+q a_j)(A_0+\alpha_2+q a_j) \big(\ell_{2j+1} \ell_{2j+2}\big)^{\varsigma_++qb_j} q^{4j+3}}, \\
a_{j+1} & \doteq A_0+q a_j , \quad \gamma_{j+1}\doteq A_1+q\gamma_j, \\
b_{j+1} &\doteq \varsigma_++qb_j , \quad \, \, \beta_{j+1}\doteq \varsigma_-+q\beta_j.
\end{align*} By using iteratively the previous relations, we may express $a_{j+1},\gamma_{j+1},b_{j+1},\beta_{j+1}$ as follows:
\begin{equation}
\begin{split}
a_{j+1} &= \left(a_0+\tfrac{A_0}{q-1}\right)q^{j+1}- \tfrac{A_0}{q-1}, \qquad \gamma_{j+1} = \left(\gamma_0+\tfrac{A_1}{q-1}\right)q^{j+1}- \tfrac{A_1}{q-1}, \\
b_{j+1} &= \tfrac{q^{j+2}-1}{q-1} \, \varsigma_+, \qquad \qquad \qquad \qquad \beta_{j+1} = \tfrac{q^{j+2}-1}{q-1}\varsigma_-.
\end{split} \label{def aj, gammaj,bj,betaj anti dS}
\end{equation} Therefore,
\begin{align*}
A_0+\alpha_{1/2}+q a_j= \alpha_{1/2}+a_{j+1} & \leqslant  \alpha_{1/2} +\left(a_0+\tfrac{A_0}{q-1}\right)q^{j+1} \\ &\leqslant  \left(\tfrac{1}{2}(b+\sqrt{b^2-4m^2}) + a_0+\tfrac{A_0}{q-1}\right)q^{j+1} \doteq M_0 q^{j+1}
\end{align*} implies $(A_0+\alpha_1+q a_j)^{-1}(A_0+\alpha_2+q a_j)^{-1}\geqslant  M_0^{-2} q^{-2(j+1)}$. Moreover, repeating similar considerations as those in the proof of Theorem \ref{Thm exponential growth}, we find that there exists a constant $M_1=M_1(\varsigma,\beta,p)>0$ such that $(\ell_{2j+1} \ell_{2j+2})^{b_{j+1}}\leqslant M_1$ holds for any $j\in\mathbb{N}$. Combining the previous estimates, we have
\begin{align*}
C_{j+1}\geqslant \frac{\mu C q(q-\tfrac{1}{2})^2}{M_0^2 M_1}q^{-6(j+1)}C_j^q \doteq G q^{-6(j+1)}C_j^q
\end{align*} for any $j\in\mathbb{N}$. From the inequality $C_j\geqslant G q^{-6j} C_{j-1}^q$, repeating analogous intermediate steps as those in the proof of Theorem \ref{Thm exponential growth}, we can find a $j_4=j_4(n,c,H,b,m^2,\mu,\varrho,\varsigma,\beta,p,R)\in\mathbb{N}$ such that for $j\geqslant j_4$ we have
\begin{align}\label{lb log Cj anti dS}
\ln C_j\geqslant q^{j}\left(\ln \big(K_1\varepsilon^q\big)-\frac{6q \ln q}{(q-1)^2}+\frac{\ln G}{q-1}\right) =q^j \ln\big(\tilde{G}\varepsilon^q\big),
\end{align} where $\tilde{G}\doteq K_1q^{-6q/(q-1)^2}G^{1/(q-1)}$.

Next, we combine \eqref{lb estimate V step j exponential growth}, \eqref{def aj, gammaj,bj,betaj anti dS} and \eqref{lb log Cj anti dS}, obtaining for $t\geqslant L\doteq \lim_{j\to\infty} L_j$ and $j\geqslant j_4$
\begin{align*}
V(t) &\geqslant \exp\left(q^j \ln\big(\tilde{G}\varepsilon^q\big)\right) \mathrm{e}^{(a_j-\gamma_j) t} (t-L)^{b_j}(1+t)^{-\beta_j} \\
& = \exp\left(q^j\left( \ln\big(\tilde{G}\varepsilon^q\big)+\left(a_0-\gamma_0+\tfrac{A_0-A_1}{q-1}\right)t+\tfrac{q\varsigma_+}{q-1}\ln(t-L)-\tfrac{q\varsigma_-}{q-1}\ln(1+t)\right)\right)\\ & \qquad \times \mathrm{e}^{\frac{A_1-A_0}{q-1} t} (t-L)^{-\tfrac{\varsigma_+}{q-1}}(1+t)^{\tfrac{\varsigma_-}{q-1}}.
\end{align*} Thus, for $t\geqslant \max\{2L,1\}$ and $j\geqslant j_4$, by $\ln(t-L)\geqslant \ln t -\ln 2$ and $-\ln(1+t)\geqslant \ln t-\ln 2$ we have
\begin{align*}
V(t) &\geqslant \exp\left(q^j\left( \ln\left(\hat{G}\varepsilon^q \mathrm{e}^{\big(a_0-\gamma_0+\tfrac{A_0-A_1}{q-1}\big)t}t^{\tfrac{\varsigma q}{q-1}}\right)\right)\right) \mathrm{e}^{\frac{A_1-A_0}{q-1} t} (t-L)^{-\tfrac{\varsigma_+}{q-1}}(1+t)^{\tfrac{\varsigma_-}{q-1}}.
\end{align*} where $\hat{G}\doteq 2^{-(\varsigma++\varsigma_-)\frac{q}{q-1}}\tilde{G}$. By \eqref{def rho crit V0} and \eqref{A0-A1}, we find 
\begin{align*}
& a_0-\gamma_0+\tfrac{A_0-A_1}{q-1}\\  & \quad  = \tfrac{q}{q-1} (\varrho-\varrho_{\mathrm{crit}}) +\tfrac{nH}{2}-\tfrac b2 -\tfrac Hp+ \tfrac{1}{q-1}\left(\tfrac 12 (b+nH)(q-1)+H(\beta+1)-nH(q-1)-\tfrac Hp\right) \\
&\quad = \tfrac{q}{q-1} (\varrho-\varrho_{\mathrm{crit}}).
\end{align*} Hence, we may rewrite the previous estimate as follows
\begin{align} \label{final lb V exponential growth anti dS nlin lb}
V(t) &\geqslant \exp\left(q^j\left( \ln\left(\hat{G}\varepsilon^q \big(\chi_{n,H,b,m^2,p,\beta,\varrho,\varsigma}(t)\big)^{\tfrac{q}{q-1}(\varrho-\varrho_{\mathrm{crit}})}\right)\right)\right) \mathrm{e}^{\frac{A_1-A_0}{q-1} t} (t-L)^{-\tfrac{\varsigma_+}{q-1}}(1+t)^{\tfrac{\varsigma_-}{q-1}}
\end{align} for $t\geqslant \max\{2L,1\}$ and $j\geqslant j_4$, where $\chi_{n,H,b,m^2,p,\beta,\varrho,\varsigma} $ is defined in \eqref{def chi function anti dS}.

By \eqref{def chi function anti dS} it follows that $\chi_{n,H,b,m^2,\beta,p,\varrho,\varsigma}$ is a strictly increasing function for $t\geqslant \bar{T}$, for some $\bar{T}=\bar{T}(n,H,b,m^2,\beta,p,\varrho,\varsigma)\geqslant 0$. Clearly, for $\varsigma\geqslant 0$ it results $\bar{T}(n,H,b,m^2,\beta,p,\varrho,\varsigma)=0$. With a slight abuse of notation, in the next lines we use the notation $\chi_{n,H,b,m^2,\beta,p,\varrho,\varsigma}^{-1}$ for the inverse function of the restriction $\chi_{n,H,b,m^2,\beta,p,\varrho,\varsigma}\big|_{[\bar{T},\infty)}$.
 
 The logarithmic factor multiplying $q^j$ in the lower bound for $V$ in \eqref{final lb V exponential growth anti dS nlin lb} is strictly positive if and only if $\hat{G}\varepsilon^q \big(\chi_{n,H,b,m^2,p,\beta,\varrho,\varsigma}(t)\big)^{\frac{q}{q-1}(\varrho-\varrho_{\mathrm{crit}})}>1$. For $t\geqslant \bar{T}$, this condition can be rewritten as
 \begin{align*}
 t> \chi_{n,H,b,m^2,p,\beta,\varrho,\varsigma}^{-1}\left((\hat{G}^{\frac{1}{q}}\varepsilon)^{-\tfrac{q-1}{\varrho-\varrho_{\mathrm{crit}}}}\right).
 \end{align*} 
 Due to $\lim_{s\to \infty} \chi_{n,H,b,m^2,p,\beta,\varrho,\varsigma}^{-1}(s)=\infty$, we can fix $\varepsilon_0=\varepsilon_0(n,c,H,b,m^2,q,\mu,\varrho,\varsigma,u_0,u_1,R)>0$  small enough that satisfies $$\chi_{n,H,b,m^2,p,\beta,\varrho,\varsigma}^{-1}\left((\hat{G}^{\frac{1}{q}}\varepsilon_0)^{-\tfrac{q-1}{\varrho-\varrho_{\mathrm{crit}}}}\right)\geqslant \max\{2L,1,\tilde{T}\}.$$ Hence, for any $\varepsilon\in(0,\varepsilon_0]$ and any $t> \chi_{n,H,b,m^2,p,\beta,\varrho,\varsigma}^{-1}\left((\hat{G}^{\frac{1}{q}}\varepsilon)^{-(q-1)/(\varrho-\varrho_{\mathrm{crit}})}\right)$ we have that $t\geqslant\{2L,1,\bar{T}\}$ and that the factor multiplying $q^j$ in \eqref{final lb V exponential growth anti dS nlin lb} is positive, therefore, letting $j\to \infty$ in \eqref{final lb V exponential growth anti dS nlin lb} we find that the lower bound for $V(t)$ is not finite. Summarizing, we showed that $V$ blows up in finite time and we obtained the upper bound estimate for the lifespan in \eqref{upper bound for the lifespan exponential growth anti dS nlin}.

\section{Final remarks} \label{Section final rmk}

In the present paper, we derived a hierarchy of blow-up results for the semilinear models \eqref{deSitter semi} and \eqref{anti deSitter semi} prescribing different levels of assumptions on the function $\Gamma(t)$.

For the semilinear Cauchy problem in de Sitter spacetime our results in Theorems \ref{Thm exponential growth}-\ref{Thm logarithmic growth} refined the results from \cite[Theorem 1.1]{Yag09} (cf. Remark \ref{Remark improvement Yagjian's result}). On the other hand, for the semilinear Cauchy problem in anti-de Sitter spacetime, to the best of our knowledge, the results from Theorems \ref{Thm anti dS lin lb}-\ref{Thm anti dS nlin lb} are completely new.

Let us make some final remarks on generalizations of the obtained results and conjectures on the models that we have treated.

\subsection{Semilinear wave equation with a summable speed of propagation}

We point out explicitly that all results we proved for the semilinear wave equation in de Sitter spacetime (namely, Theorems \ref{Thm exponential growth}, \ref{Thm polynomial growth} and \ref{Thm logarithmic growth}) can be naturally extended to the following semilinear Cauchy problem 
\begin{align*}
\begin{cases}
\partial_t^2 u- a^2(t)\Delta u+ b\partial_t u +m^2 u=f(t,u), & x\in \mathbb{R}^n, \ t\in (0,T), \\
u(0,x)= \varepsilon u_0(x), & x\in \mathbb{R}^n, \\
\partial_t u(0,x)= \varepsilon u_1(x), & x\in \mathbb{R}^n,
\end{cases}
\end{align*}
where $a\in L^1([0,\infty))$ is a nonnegative function and $f(t,u)$ is given by \eqref{nonlinearity}. This is straightforward consequence of the fact that $A(t)=\int_0^t a(\tau) \mathrm{d}\tau$ is a bounded function. 

Nevertheless, in the previous sections we consider just the case $a=a_{\mathrm{dS}}$, since we expect that our results cannot be improved by working with the spatial average of a local solution as we explained in Remark \ref{Remark V0 for de Sitter}.

\subsection{Dominant mass case}

In the case with dominant mass ($b^2<4m^2$), our approach is unfruitful due to the conjugate complex roots $\alpha_{1/2}$ in \eqref{decomposition 2nd order operator}. For the average of a local solution $u$ to \eqref{deSitter semi}, we obtain the representation
\begin{align*}
U(t)& = \varepsilon \, y_0(t;b,m^2) \int_{\mathbb{R}^n}u_0(x)\, \mathrm{d}x +\varepsilon \, y_1(t;b,m^2) \int_{\mathbb{R}^n}u_1(x)\, \mathrm{d}x \\ & \qquad + \int_0^t y_1(t-\tau;b,m^2) \, \Gamma(\tau) \left(\int_{\mathbb{R}^n}|u(\tau,x)|^p\mathrm{d}x \right)^{\beta+1} \mathrm{d}\tau,
\end{align*} where 
\begin{align*}
y_0(t;b,m^2) &\doteq \mathrm{e}^{-\frac{b}{2}t}  \cos\left(\sqrt{m^2-\frac{b^2}{4}} t\right)+\frac{b}{2}\frac{\mathrm{e}^{-\frac{b}{2}t}  \sin\left(\sqrt{m^2-\frac{b^2}{4}} t\right)}{\sqrt{m^2-\frac{b^2}{4}}} ,\\
y_1(t;b,m^2) &\doteq \frac{\mathrm{e}^{-\frac{b}{2}t}  \sin\left(\sqrt{m^2-\frac{b^2}{4}} t\right)}{\sqrt{m^2-\frac{b^2}{4}}}.
\end{align*} Analogously, for a local solution $v$ to \eqref{anti deSitter semi}. We see that the damped oscillations from the time factors $y_0,y_1$ prevent us to work with a nonnegative $U$ making impossible to establish an iteration frame and, consequently, a sequence of lower bound estimates for $U$.

\subsection{Critical exponent for anti-de Sitter with a nonlocal nonlinearity}
In the main results from Subsection \ref{Subsection Main results} we analyzed how prescribing certain conditions on $\Gamma$ it is possible to prove blow-up results for local solutions to \eqref{deSitter semi} and \eqref{anti deSitter semi}. Let us now change our perspective in the following sense: is it possible to prove blow-up results either for \eqref{deSitter semi} or for \eqref{anti deSitter semi} without requiring additional exponential/polynomial growth for the nonlinear term though the factor $\Gamma$? As we explained in the introduction, for the model \eqref{deSitter semi} this is not possible unless $m^2=0$, and this was actually one of the reasons for us to consider a nonlinear term given by \eqref{nonlinearity}. Similarly, for \eqref{anti deSitter semi} when \eqref{1st lb for V from lin part is dominant} holds, due to the fact that $\varrho_{\mathrm{crit}}>0$, Theorems \ref{Thm anti dS lin lb}, \ref{Thm anti dS lin lb poly growth} and \ref{Thm anti dS lin lb log growth} do not provide a blow-up result unless $\varrho>0$. On the other hand, when \eqref{1st lb for V from V0 is dominant} holds it is possible to find some $p>1$ such that $\varrho_{\mathrm{crit}}\leqslant 0$. Indeed, from straightforward computations we see that the condition $\varrho_{\mathrm{crit}}\leqslant 0$ can be rewritten as the following quadratic equation for $p$
\begin{align}\label{quadratic eq p}
 (b+nH)(\beta+1)(p-1)^2+\left[2(b+H)\beta +(b+nH)+H\right] (p-1)- 
((n-2)H-b)\beta\leqslant 0.
\end{align} In the nonlocal case $\beta>0$, the term $-((n-2)H-b)\beta$ in \eqref{quadratic eq p} can be negative, more precisely this happens for $n > 2+ \frac{b}{H}$. Therefore,  for $n > 2+ \frac{b}{H}$, thanks to Descartes' rule of signs, there exists $p_0=p_0(n,H,b,\beta)>1$ such that for $1<p\leqslant p_0(n,H,b,\beta)$ the quadratic equation \eqref{quadratic eq p} is satisfied. Since $2+\frac{b}{H}\geqslant N_0$, where $N_0$ is defined in \eqref{def N0 dim}, in particular we are in the case in which \eqref{1st lb for V from V0 is dominant} is always fulfilled for any $p>1$, hence, Theorem \ref{Thm anti dS nlin lb} provides a blow-up result for $1<p< p_0(n,H,b,\beta)$ for local in time solutions to \eqref{anti deSitter semi} when $n> 2+\frac{b}{H}$. 

In other words, when $n> 2+\frac{b}{H}$ and $\beta>0$, considering $\Gamma(t)=1$, that is, when the nonlinearity in \eqref{anti deSitter semi} is given by 
\begin{align*}
f(v)= \left(\int_{\mathbb{R}^n}|v(t,y)|^p\mathrm{d}y\right)^\beta |v|^p,
\end{align*} we found as candidate to be the critical exponent for the semilinear Cauchy problem \eqref{anti deSitter semi} the positive root $p_0(n,H,b,\beta)$ of the quadratic equation
 \begin{align*}
 (b+nH)(\beta+1)(p-1)^2+\left[2(b+H)\beta +(b+nH)+H\right] (p-1)- 
((n-2)H-b)\beta = 0.
\end{align*} In \cite{PalTak22crit} we will prove a blow-up result for local solutions to \eqref{anti deSitter semi} for $\varrho=\varrho_{\mathrm{crit}}$ when \eqref{1st lb for V from V0 is dominant} holds and $\varrho_{\mathrm{crit}}$ is given by \eqref{def rho crit V0}. In particular, for $n> 2+\frac{b}{H}$ and $\beta>0$ this result will show the blow-up of local solutions even in the threshold case $p=p_0(n,H,b,\beta)$.


\section*{Acknowledgments}

A. Palmieri is supported by the\emph{ Japan Society for the Promotion of Science} (JSPS) – JSPS Postdoctoral Fellowship for Research in Japan (Short-term) (PE20003) – and is member of the \emph{Gruppo Nazionale per l’Analisi Matematica, la Probabilit\`{a} e le loro Applicazioni} (GNAMPA) of the \emph{Instituto Nazionale di Alta Matematica} (INdAM). H. Takamura is partially supported by the Grant-in-Aid for Scientific Research (B) (No.18H01132), \emph{Japan Society for the Promotion of Science}.

\end{document}